\documentclass{amsart}
\usepackage{graphicx}

\newtheorem{theorem}{Theorem}[section]
\newtheorem{lemma}[theorem]{Lemma}

\theoremstyle{definition}

\newtheorem{example}[theorem]{Example}
\newtheorem{xca}[theorem]{Exercise}
\newtheorem{claim}[theorem]{Claim}

\theoremstyle{remark}
\newtheorem{remark}[theorem]{Remark}

\numberwithin{equation}{section}

\begin{document}

\title{Morse position of knots and closed incompressible surfaces}

\author{Makoto Ozawa}
\address{Department of Natural Sciences, Faculty of Arts and Sciences, Komazawa University,
 1-23-1 Komazawa, Setagaya-ku, Tokyo, 154-8525, Japan}
\email{w3c@komazawa-u.ac.jp}





\begin{abstract}
In this paper, we study on knots and closed incompressible surfaces in the 3-sphere via Morse functions.
We show that both of knots and closed incompressible surfaces can be isotoped into a "related Morse position" simultaneously.
As an application, we have following results.
\begin{itemize}
\item Smallness of Montesinos tangles with length two and Kinoshita's theta curve
\item Classification of closed incompressible and meridionally incompressible surfaces in 2-bridge theta-curve and handcuff graph complements and the complements of links which admit Hopf tangle decompositions
\end{itemize}
\end{abstract}

\maketitle

\section{Introduction}
Morse theory had a tremendous impact on low-dimensional topology.
Bridge positions of knots that developed by Schubert (\cite{S1}) and Heegaard splittings of 3-manifolds that developed by Heegaard (\cite{H}) are independent of Morse theory, but those can be explained integrally by Morse theory.


A bridge position of a knot is defined by using a Morse function from the 3-sphere to the 1-dimensional Euclidean space with two critical points so that any maximal critical point of the knot is situated above any minimal critical point.
This general idea was expanded to a notion of thin position by Gabai, and he settled the Property R conjecture (\cite{G}).
Thin position also brought a celebrated solution of the knot complement conjecture by Gordon and Luecke (\cite{GL})

In this paper, we show that we can isotope a closed incompressible surface to a related Morse position in the complement of a knot placed in a bridge or thin position.
Then as an applications, we show smallness of Montesinos tangles with length two and Kinoshita's theta curve, and classification of closed incompressible and meridionally incompressible surfaces in 2-bridge theta-curve and handcuff graph complements and the complements of links which admit Hopf tangle decompositions.


\section{Result}

Let $K$ be a knot in the 3-sphere $S^3$ and $h:S^3\to \Bbb{R}$ be the composition of the inclusion $S^3\subset \Bbb{R}^4$ and the height function $\Bbb{R}^4\to \Bbb{R}$ defined by $(x,y,z,w)\mapsto w$.
We say that $K$ is in a {\em Morse position} with respect to $h$ if the restriction of $h$ to $K$ is a Morse function.
Let $a_0,\ldots,a_n$ be the critical points of $K$ which are labeled so that the corresponding critical values $t_i=h(a_i)$ satisfy $t_{i-1}<t_i$ for all $i$.
For a regular value $s_i\in(t_{i-1},t_i)$, $P_i=h^{-1}(s_i)$ is called a {\em level sphere} between $a_{i-1}$ and $a_i$.
The {\em width} of $K$ is $w(K)=\sum_{i=1}^{n} |P_i\cap K|$.
A knot $K$ is in {\em thin position} if $w(K)$ is minimal up to isotopy of $K$.
A level sphere $P_i$ between $a_{i-1}$ and $a_i$ is said to be {\em thin} if $a_{i-1}$ is a local maximum and $a_i$ a local minimum.
Similarly, a level sphere $P_i$ between $a_{i-1}$ and $a_i$ is said to be {\em thick} if $a_{i-1}$ is a local minimum and $a_i$ a local maximum.
A knot $K$ is in {\em bridge position} if there is only one thick level sphere $P_{\frac{n+1}{2}}$.
Then, the {\em bridge number} of $K$ is $b(K)=|P_{\frac{n+1}{2}}\cap K|/2$.
A knot $K$ is in {\em minimal bridge position} if $b(K)$ is minimal up to isotopy of $K$.
For each thick level sphere $P_i$ between $a_{i-1}$ and $a_i$, we define the {\em thick region} of $K$ as $h^{-1}([t_{i-1}+\epsilon, t_i-\epsilon])$ for a sufficiently small positive real number $\epsilon$.
Then, no thick region contains a maximum or minimum of $K$, and the rest of all thick regions in $S^3$ contain all maxima and minima of $K$, each of which is called a {\em thin region} of $K$.
Any thick region contains only one thick level sphere, and any thin region, except for the top 3-ball and the bottom 3-ball, contains only one thin level sphere.
Let $M_{thin}$ be the disjoint union of all thin regions of $K$, and $M_{thick}$ the disjoint union of all thick regions of $K$.

Let $F$ be a closed surface in $S^3$ disjoint from $K$ or intersecting $K$ transversely.
We say that $F$ is in a {\em Morse position} with respect to $h$ if the restriction of $h$ to $F$ is a Morse function.
The critical points of $K$ and $F$ are assumed to be mutually disjoint.
Moreover, we say that $F$ is in a {\em Morse position related} to $K$ if the following conditions are satisfied.
\begin{enumerate}
\item[(1)] $F$ is in a Morse position with respect to $h$,
\item[(2)] $F$ and $K$ are in general position and $F\cap K$ is contained in $M_{thick}$,
\item[(3)] all maxima and minima of $F$ are contained in $M_{thin}$, and all saddles of $F$ are contained in $M_{thick}$.
\end{enumerate}

These three conditions are easily satisfied for any closed surface $F$.
We prepare definitions to describe the next condition (4).

Let $M$ be a 3-manifold, $T$ a 1-manifold properly embedded in $M$, and $F$ a surface properly embedded in $M$ such that $F$ is disjoint from $T$ or intersects $T$ transversely in the interior of $F$.
An {\em isotopy} $\phi_t$ of $F$ {\em in} $(M,T)$ should be assumed that $\phi_t(F)$ and $T$ are in general position for all $t\in[0,1]$.
We recall the definitions about incompressible surfaces and meridionally incompressible surfaces.

First, suppose that there exists a disk $D$ in $M-T$ such that $D\cap F=\partial D$ and $\partial D$ is essential in $F-T$.
Then by cutting $F$ along $\partial D$ and pasting two parallel copies of $D$, we have a new surface $F'$.
We say that $F'$ is obtained from $F$ by {\em compressing $F$ along $D$}, and $D$ is called a {\em compressing disk} for $F$.
A surface $F$ is {\em incompressible} in $(M,T)$ if for any disk $D$ in $M-T$ with $D\cap F=\partial D$, there exists a disk $D'$ in $F-T$ such that $\partial D'=\partial D$.
In the case that $F$ is a disk disjoint from $T$, $F$ is said to be {\em incompressible} in $(M,T)$ if for any disk $D$ in $\partial M-\partial T$ with $\partial D=\partial F$, $F\cup D$ does not bound a 3-ball in $M-T$.
In the case that $F$ is a 2-sphere disjoint from $T$, $F$ is said to be {\em incompressible} in $(M,T)$ if $F$ does not bound a 3-ball in $M-T$.

Next, suppose that there exists a disk $D$ in $M$ such that $D\cap F=\partial D$ and $D$ intersects $T$ in one point of the interior of $D$.
Then by cutting $F$ along $\partial D$ and pasting two parallel copies of $D$, we have a new surface $F'$.
If $F$ is incompressible in $(M,T)$, then $F'$ is also incompressible in $(M,T)$.
We say that $F$ is {\em meridionally incompressible} in $(M,T)$ if for any disk $D$ in $M$ such that $D\cap F=\partial D$ and $D$ intersects $T$ in one point of the interior of $D$, there exists a disk $D'$ in $F$ such that $\partial D'=\partial D$ and $D\cup D'$ bounds a trivial $($3-ball, arc$)$-pair in $(M,T)$.
In the case that $F$ is a disk intersecting $T$ in one point in the interior of $F$, we say that $F$ is {\em meridionally incompressible} in $(M,T)$ if for any disk $D$ in $\partial M$ with $\partial D=\partial F$, $F\cup D$ does not bound a trivial $($3-ball, arc$)$-pair in $($M,T$)$.
In the case that $F$ is a 2-sphere intersecting $T$ in two points, we say that $F$ is {\em meridionally incompressible} in $(M,T)$ if $F$ does not bound a trivial $($3-ball, arc$)$-pair in $(M,T)$.
We note that this definition is slightly different from a usual one as a matter of convenience.

Suppose that $K\subset S^3$ is a knot in a Morse position with respect to $h$ and $F\subset S^3$ is a closed surface in a Morse position related to $K$.
Furthermore, we say that $F$ is in an {\em essential Morse position related} to $K$ if the following condition is satisfied.
\begin{enumerate}
	\item [(4)] each component of $F\cap M_{thin}$ and $F\cap M_{thick}$ is incompressible and meridionally incompressible in $(M_{thin},K\cap M_{thin})$ and $(M_{thick},K\cap M_{thick})$ respectively.
\end{enumerate}


The following main theorem has a philosophy that makes thin regions $M_{thin}$ simple, but thick regions $M_{thick}$ essential.

\bigskip
\begin{theorem}
\label{main}
{\it Let $K\subset S^3$ be a knot in a Morse position with respect to $h$ and $F\subset S^3$ be a closed surface disjoint from $K$ or intersecting $K$ transversely such that $F$ is incompressible and meridionally incompressible in $(S^3,K)$.
Then $F$ can be isotoped so that;
\begin{enumerate}
	\item $F$ is a thin level sphere, or
	\item $F$ is in an essential Morse position related to $K$.
\end{enumerate}
}
\end{theorem}


\bigskip
\begin{example}
Let $K$ be a composite knot of two trefoils, and $F$ be the decomposing sphere for $K$.
If we put $K$ in a bridge position, then by Theorem \ref{main}, $F$ can be isotoped so that $F$ is in an essential Morse position related to $K$.
See Figure \ref{type 2}.

\begin{figure}[htbp]
	\begin{center}
		\includegraphics[trim=0mm 0mm 0mm 0mm, width=.5\linewidth]{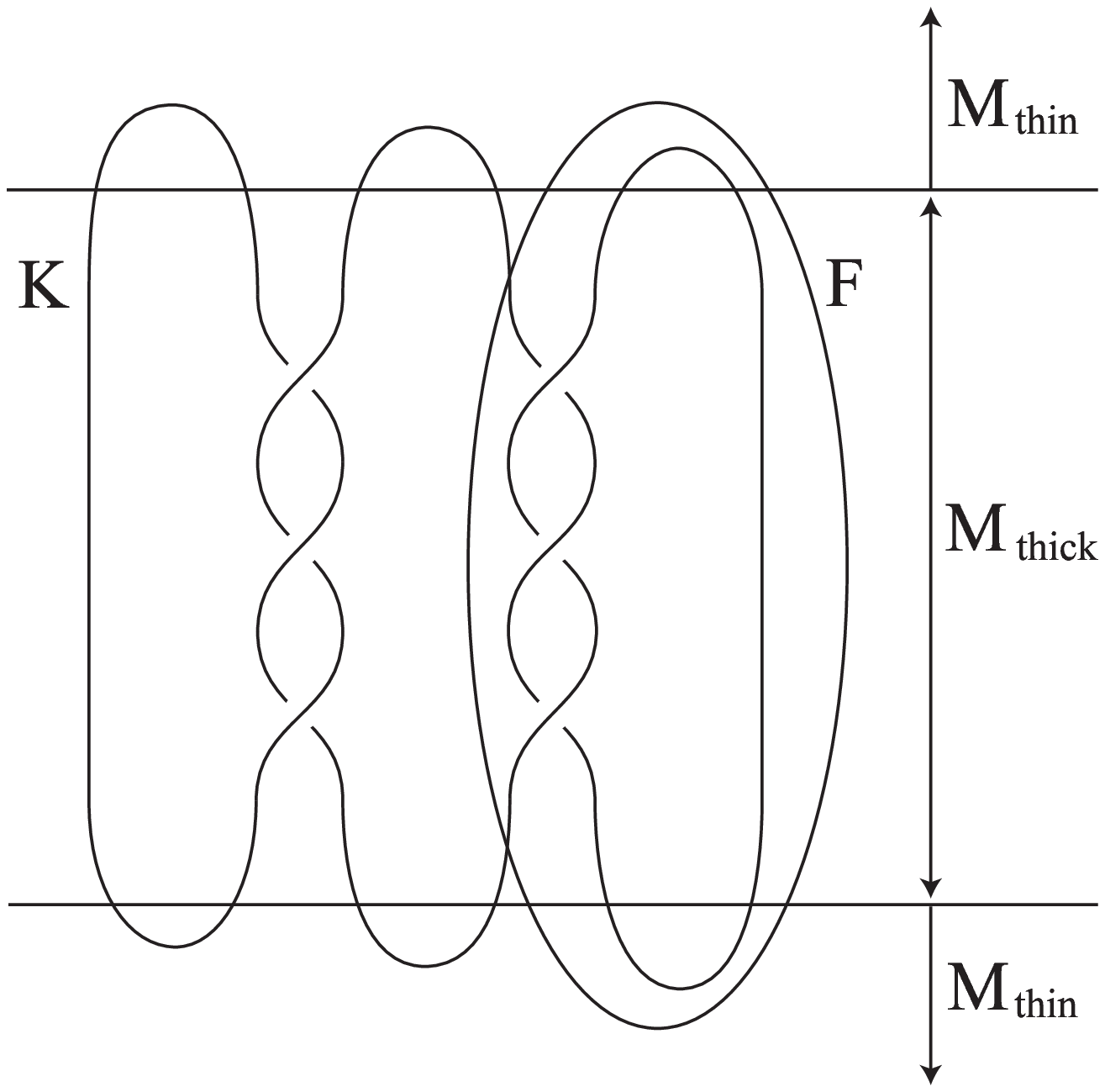}
	\end{center}
	\caption{Simple example of type 2 in Theorem \ref{main}}
	\label{type 2}
\end{figure}

If we put $K$ in a thin position, then by Theorem \ref{main}, $F$ can be isotoped so that $F$ is in a thin level sphere.
See Figure \ref{type 1}.

\begin{figure}[htbp]
	\begin{center}
		\includegraphics[trim=0mm 0mm 0mm 0mm, width=.5\linewidth]{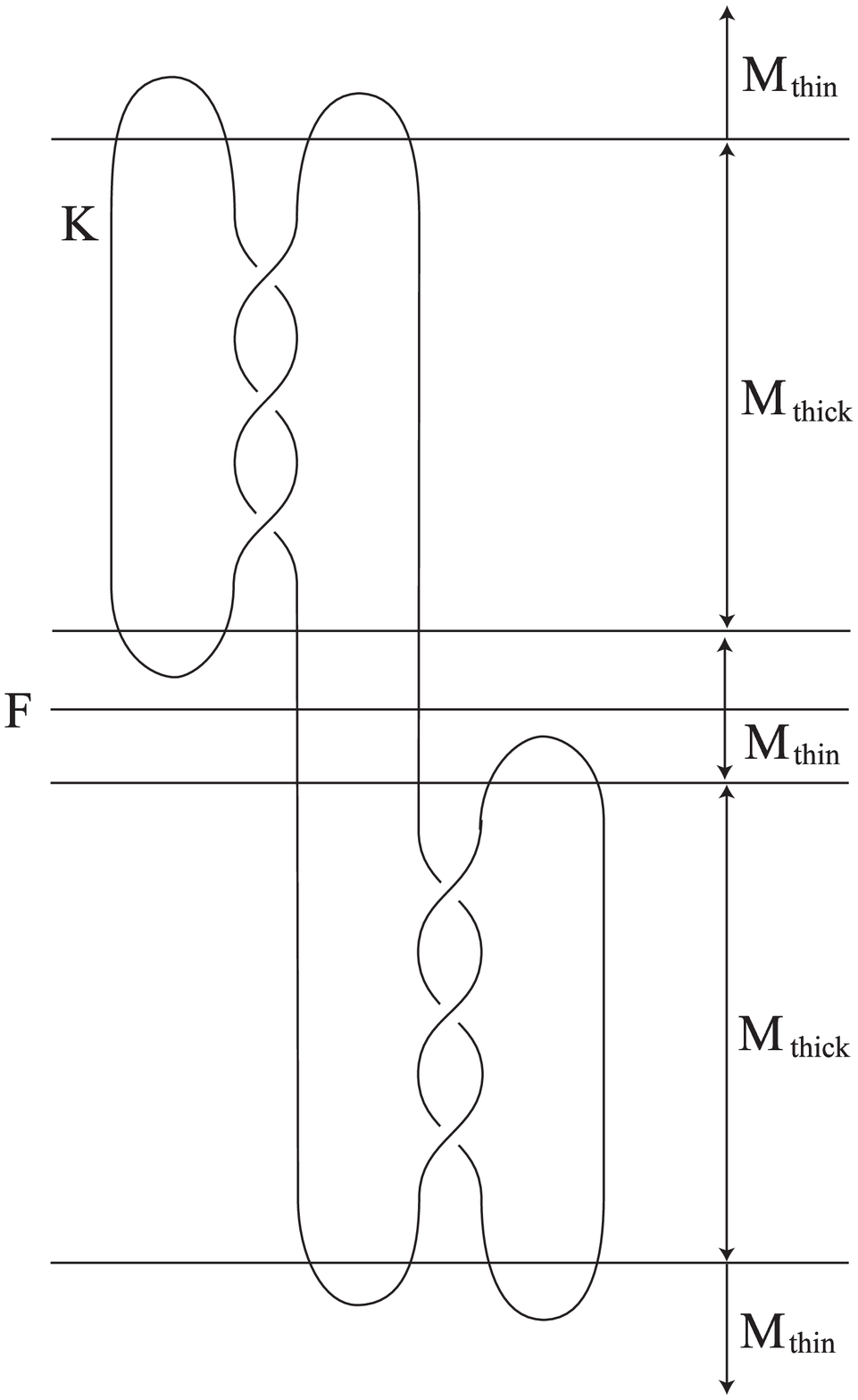}
	\end{center}
	\caption{Simple example of type 1 in Theorem \ref{main}}
	\label{type 1}
\end{figure}
\end{example}

A 3-manifold $M$ is said to be {\em small} if any closed incompressible surface properly embedded in $M$ is $\partial$-parallel.
Analogously, a knot $K\subset S^3$ is {\em small} if the exterior $E(K)=S^3-{\rm int} N(K)$ is small, a tangle $(B,T)$ is {\em small} if $B-{\rm int} N(T)$ is small, and a graph $G\subset S^3$ is {\em small} if $S^3-{\rm int} N(G)$ is small.

\begin{figure}[htbp]
	\begin{center}
		\includegraphics[trim=0mm 0mm 0mm 0mm, width=.5\linewidth]{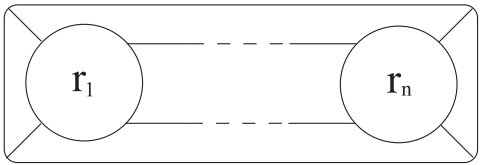}
	\end{center}
	\caption{Montesinos tangle $T(r_1,\ldots,r_n)$}
	\label{Montesinos tangle}
\end{figure}

A {\em Montesinos tangle} $T(r_1,\ldots,r_n)$ is defined as a series of $n$ rational tangles of slope $r_i$ $(i=1,\ldots,n)$.
See Figure \ref{Montesinos tangle}.
We note that if $r_i=\infty$ for some $i$, then $T(r_1,\ldots,r_n)$ is splittable by a properly embedded disk,
and that if $r_i$ is an integer, then the length $n$ of $T(r_1,\ldots,r_n)$ can be taken short.
For this reason, it is reasonable to assume that $r_i\ne \infty$ and $r_i$ is not an integer for all $i$.

\bigskip
\begin{theorem}\label{Montesinos}
{\it A Montesinos tangle $T(r_1,r_2)$ is small if $r_i\ne \infty$ for $i=1,2$.}
\end{theorem}

\bigskip
\begin{remark}
The referee suggested the possibility to show that a tangle is small by showing the double branched cover is small.
Indeed, it seems true and Theorem \ref{Montesinos} follows the smallness of the double branched cover over a Montesinos tangle $T(r_1,r_2)$.
Theorem \ref{Montesinos} gives an another direct proof without branched covering spaces or orbifolds (\cite[Theorem 1]{O}).
\end{remark}

\bigskip
\begin{example}
As a corollary of Theorem \ref{Montesinos}, we can show that Kinoshita's theta curve is small.
Kinoshita's theta curve is known as a most popular theta curve, which is normally illustrated as Figure \ref{Kinoshita1}.
It is known to be minimally knotted and hyperbolic.
Therefore, the exterior is $\partial$-irreducible.

\begin{figure}[htbp]
	\begin{center}
		\includegraphics[trim=0mm 0mm 0mm 0mm, width=.4\linewidth]{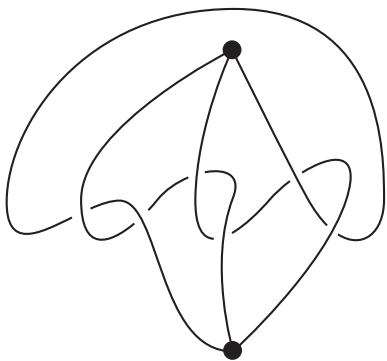}
	\end{center}
	\caption{standard form of Kinoshita's theta curve}
	\label{Kinoshita1}
\end{figure}

We shorten an edge and obtain an another form of Kinoshita's theta curve as in Figure \ref{Kinoshita2}.
Then, the 2-string tangle on the left-hand side is a Montesinos tangle with length two and the right-hand side is a trivial $H$-shaped graph tangle.
Hence, the exterior of Kinoshita's theta curve is homeomorphic to the exterior of the Montesinos tangle.
This fact was also observed by Wu (\cite[Example 3]{Wu}).
Therefore, Kinoshita's theta curve is small by Theorem \ref{Montesinos}.

\begin{figure}[htbp]
	\begin{center}
		\includegraphics[trim=0mm 0mm 0mm 0mm, width=.4\linewidth]{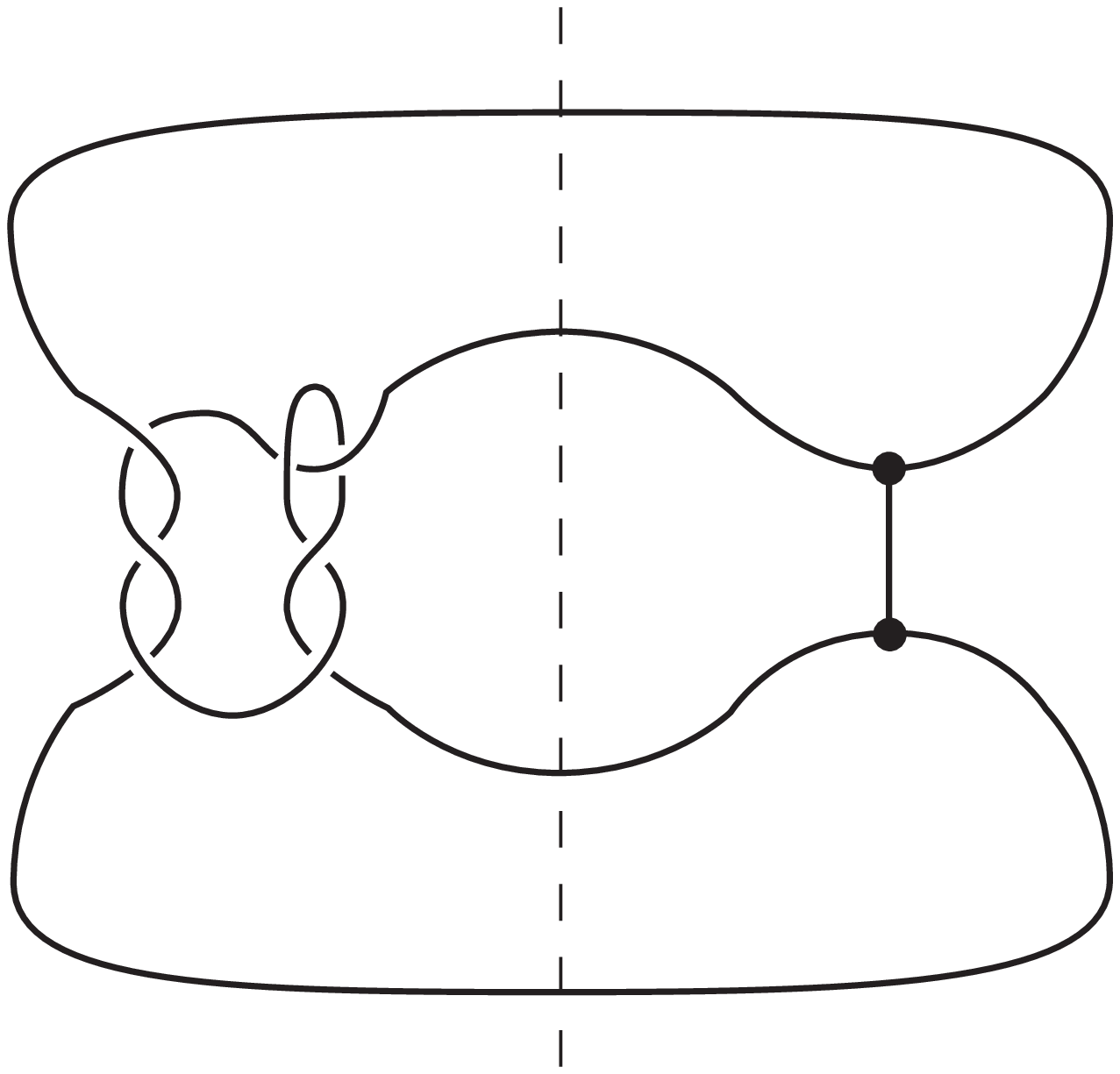}
	\end{center}
	\caption{another form of Kinoshita's theta curve}
	\label{Kinoshita2}
\end{figure}
\end{example}

Let $G$ be a theta curve or handcuff graph embedded in $S^3$.
We say that $G$ is in a {\em Morse position} with respect to $h$ if the restriction of $h$ to each edge of $G$ is a Morse function.
In this paper, we assume that near each vertex $v$ of $G$, three edges either go up or go down from $v$.
We say that a theta curve or handcuff graph $G$ in a Morse position is {\em 2-bridge} if there exists a level 2-sphere $P$ which intersects $G$ in 5 points other than vertices and critical points, and edges of $G$ have exactly two critical points.

Let $F\subset S^3$ be a closed surface disjoint from $G$ or intersecting edges of $G$ transversely.
Concerning the definition for $F$ to be meridionally incompressible, we say that $F$ is {\em meridionally incompressible} in $(S^3,G)$ if $F$ is not parallel to $\partial N(v)$ for a vertex $v$ of $G$, together with the definition for knot case.

The next theorem asserts there is no incompressible and meridionally incompressible surface of positive genus in the complement of 2-bridge theta curves.

\begin{figure}[htbp]
	\begin{center}
		\includegraphics[trim=0mm 0mm 0mm 0mm, width=.7\linewidth]{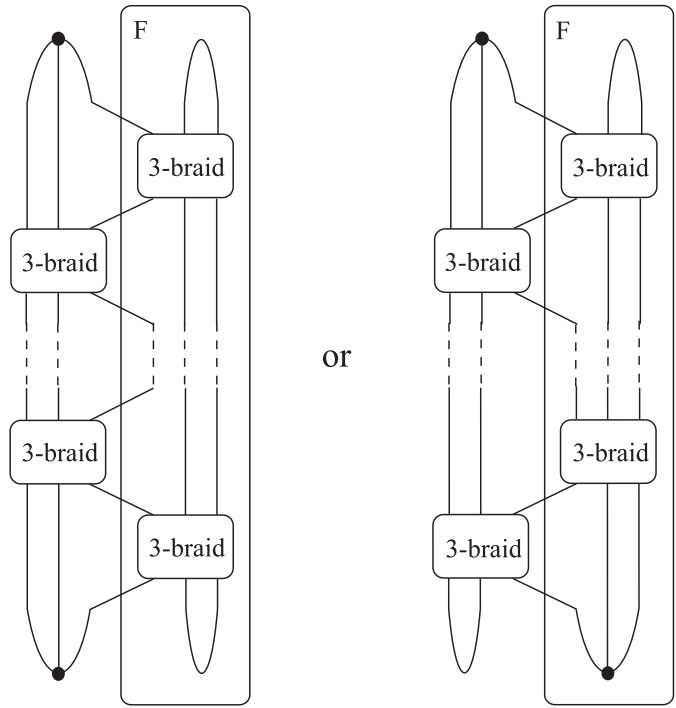}
	\end{center}
	\caption{incompressible and meridionally incompressible surface in 2-bridge theta curve or handcuff graph complements}
	\label{theta-figure}
\end{figure}

\bigskip
\begin{theorem}\label{theta}
{\it Let $G$ be a 2-bridge theta curve or handcuff graph in $S^3$, and $F$ be an incompressible and meridionally incompressible surface.
Then $F$ is isotopic to a sphere with one maximum and one minimum illustrated in Figure \ref{theta-figure}.}
\end{theorem}

\begin{figure}[htbp]
	\begin{center}
		\includegraphics[trim=0mm 0mm 0mm 0mm, width=.25\linewidth]{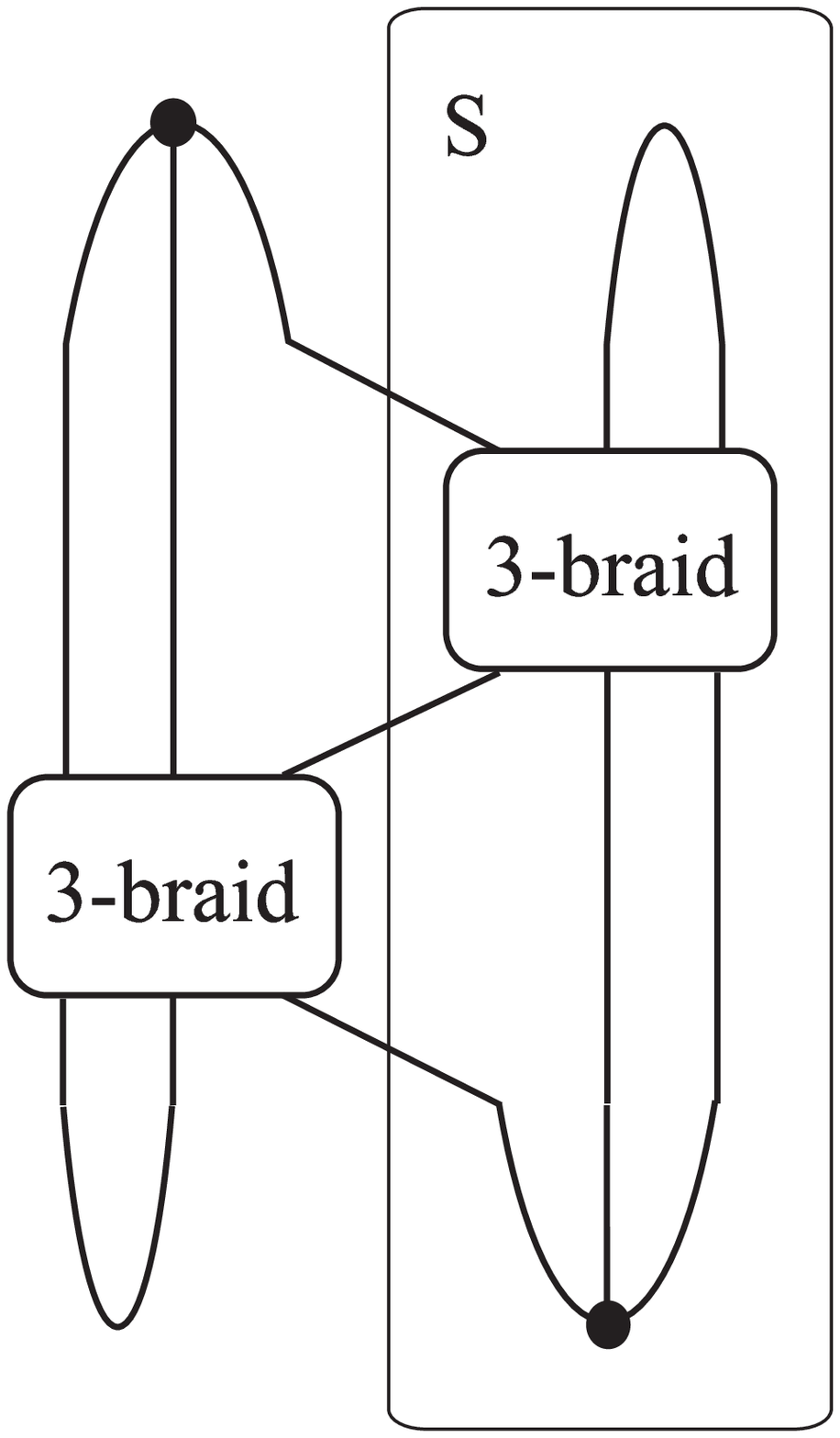}
	\end{center}
	\caption{composite 2-bridge theta curve and the decomposing sphere}
	\label{composite}
\end{figure}

As a corollary of Theorem \ref{theta}, we obtain the Motohashi's result (\cite{M}) about composite 2-bridge theta curves.
A theta curve $G$ is {\em composite} if there exists a 2-sphere which intersects each edge of $G$ in one point and decompose $G$ into two non-trivial theta curve.
By Theorem \ref{theta}, the decomposing sphere is isotopic to a sphere with one maximum and minimum as Figure \ref{composite}.
Hence, $G$ is decomposed into two rational theta curves, that is, a theta curve obtained by pasting a rational tangle and a trivial $H$-shaped graph tangle.

\begin{figure}[htbp]
	\begin{center}
		\includegraphics[trim=0mm 0mm 0mm 0mm, width=.35\linewidth]{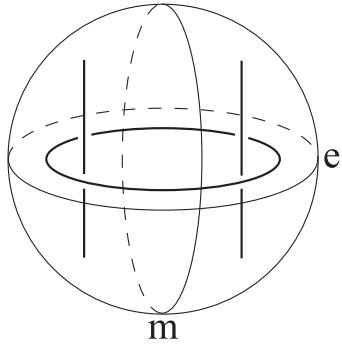}
	\end{center}
	\caption{Hopf tangle $(B,T)$}
	\label{Hopf_fig}
\end{figure}

Next, we consider closed incompressible surfaces in the complements of links that admit Hopf tangle decompositions.
The {\em Hopf tangle} $(B,T)$ is a 2-string trivial tangle $(B,t_1\cup t_2)$ with a trivial loop $C$ as illustrated as Figure \ref{Hopf_fig}.
The {\em equator} of the Hopf tangle $(B,T)$ is defined as a loop $e$ in $\partial B$ which cobounds an annulus in $B-(t_1\cup t_2)$ with the trivial loop $C$, and the {\em meridian} is defined as a loop $m$ in $\partial B$ which bounds a disk $D$ in $B-(t_1\cup t_2)$ separating two strings of $T$.

Let $\mathcal{L}(x_1,x_2)$ be the set of links which are obtained by gluing two Hopf tangles $(B_1,T_1)$ and $(B_2,T_2)$ so that $x_1$ coincides with $x_2$, where $x_i$ represents either the equator $e_i$ of $(B_i,T_i)$ or the meridian $m_i$ of $(B_i,T_i)$ for $i=1,2$.

We note that if a link $L$ in $S^3$ admits a Hopf tangle decomposition $(S^3,L)=(B_1,T_1)\cup (B_2,T_2)$, then the decomposing sphere $S=\partial B_1=\partial B_2$ is incompressible and meridionally incompressible in $(S^3,L)$.


\bigskip
\begin{remark}\label{class_remark}
We remark that $\mathcal{L}(m_1,m_2)\cap (\mathcal{L}(m_1,e_2)\cup \mathcal{L}(e_1,m_2))=\emptyset$ and $(\mathcal{L}(m_1,e_2)\cup \mathcal{L}(e_1,m_2))\cap \mathcal{L}(e_1,e_2)=\emptyset$, and that $\mathcal{L}(m_1,m_2)\cap \mathcal{L}(e_1,e_2)$ consists of one link which is obtained by the orientation reversing identity map between two Hopf tangles, namely, the Borromean rings.
\end{remark}

\bigskip
\begin{theorem}\label{Hopf}
{\it Let $L$ be a link in $S^3$ which admits a Hopf tangle decomposition $(S^3,L)=(B_1,T_1)\cup (B_2,T_2)$.
Then, 

\begin{enumerate}
	\item There exists an incompressible and meridionally incompressible closed surface in $(S^3,L)$ different from the decomposing sphere if and only if $L\in \mathcal{L}(m_1,m_2)\cup \mathcal{L}(m_1,e_2)\cup \mathcal{L}(e_1,m_2)\cup \mathcal{L}(e_1,e_2)$.
	\item $L$ admits a unique Hopf tangle decomposition if and only if $L\not\in\mathcal{L}(m_1,m_2)\cup \mathcal{L}(m_1,e_2)\cup \mathcal{L}(e_1,m_2)$.
	If $L\in \mathcal{L}(m_1,m_2)\cup \mathcal{L}(m_1,e_2)\cup \mathcal{L}(e_1,m_2)$ except the Borromean rings, then $L$ admits exactly two Hopf tangle decompositions.
	The Borromean rings admits exactly three Hopf tangle decompositions.
	\item $L$ is small if and only if $L\not\in \mathcal{L}(e_1,e_2)$.
	If $L\in \mathcal{L}(e_1,e_2)$, then there exists an essential torus in $S^3-\rm{int}N(L)$.
\end{enumerate}
}
\end{theorem}

\begin{figure}[htbp]
	\begin{center}
	\begin{tabular}{ccc}
		\includegraphics[trim=0mm 0mm 0mm 0mm, width=.25\linewidth]{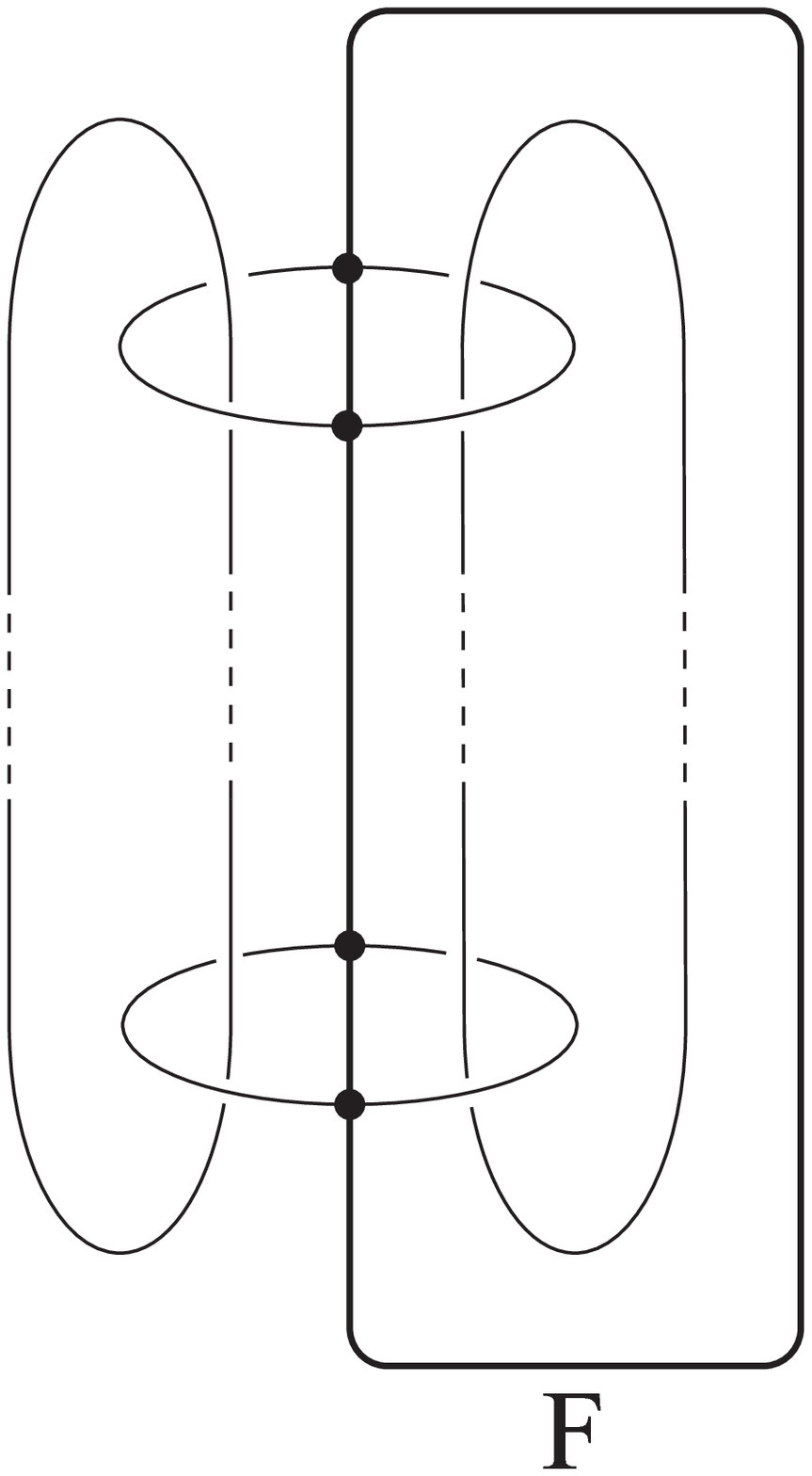} &
		\includegraphics[trim=0mm 0mm 0mm 0mm, width=.25\linewidth]{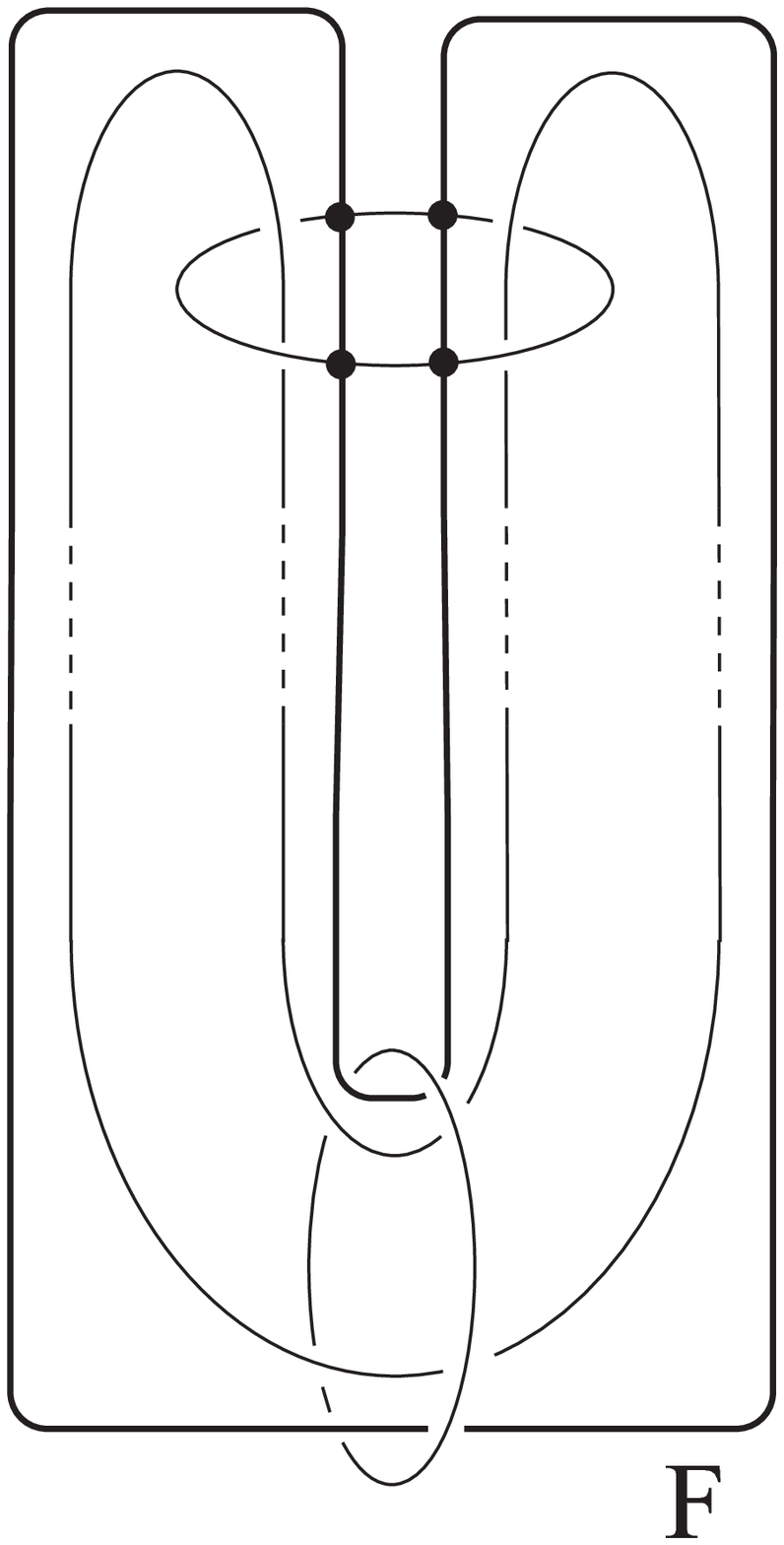} &
		\includegraphics[trim=0mm 0mm 0mm 0mm, width=.25\linewidth]{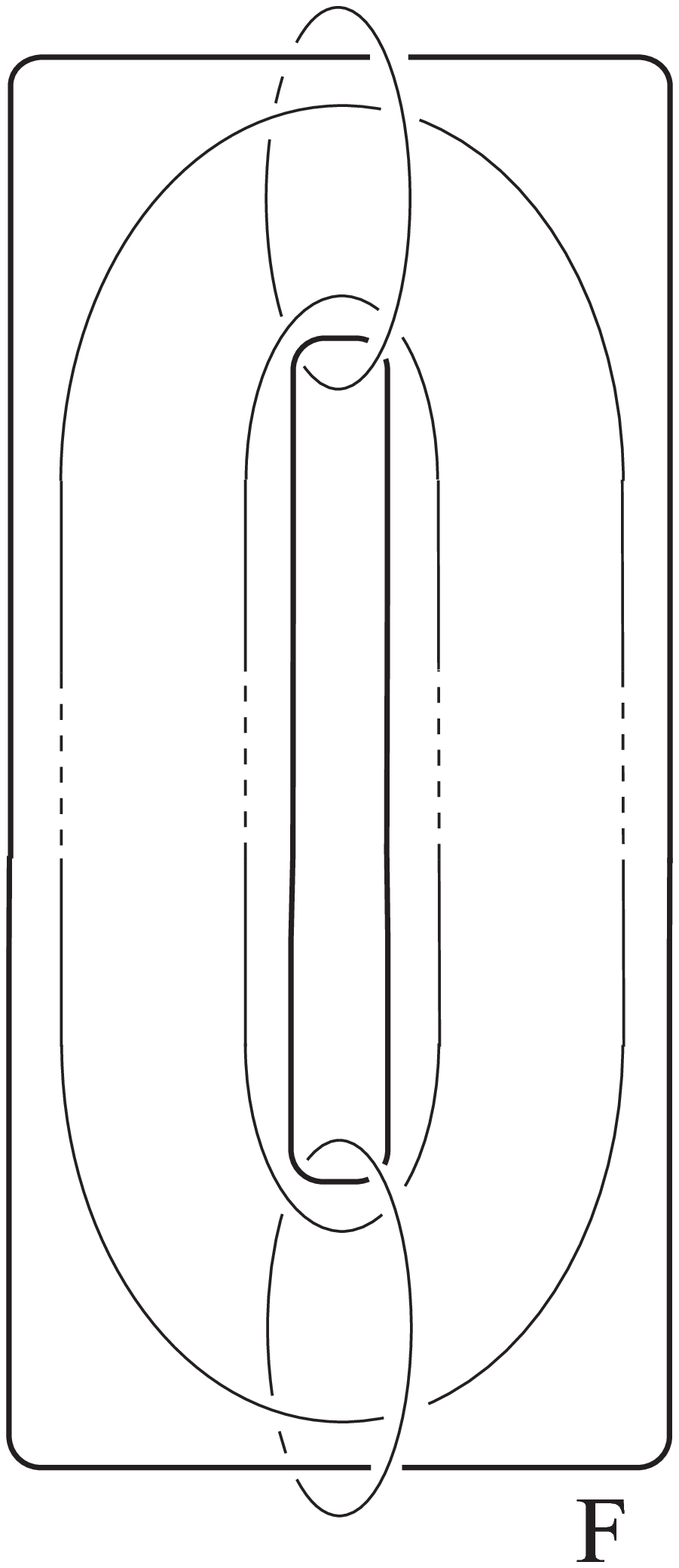} \\
		$\mathcal{L}(m_1,m_2)$ & $\mathcal{L}(m_1,e_2)\cup \mathcal{L}(e_1,m_2)$ & $\mathcal{L}(e_1,e_2)$
	\end{tabular}
	\end{center}
	\caption{three types of link classes and closed incompressible and meridionally incompressible surfaces}
	\label{class}
\end{figure}

See Figure \ref{class} for the types of link classes and associated closed incompressible and meridionally incompressible surfaces, where dotted lines represent monotone strings (which may have some twists).


\begin{figure}[htbp]
	\begin{center}
	\begin{tabular}{cc}
		\includegraphics[trim=0mm 0mm 0mm 0mm, width=.4\linewidth]{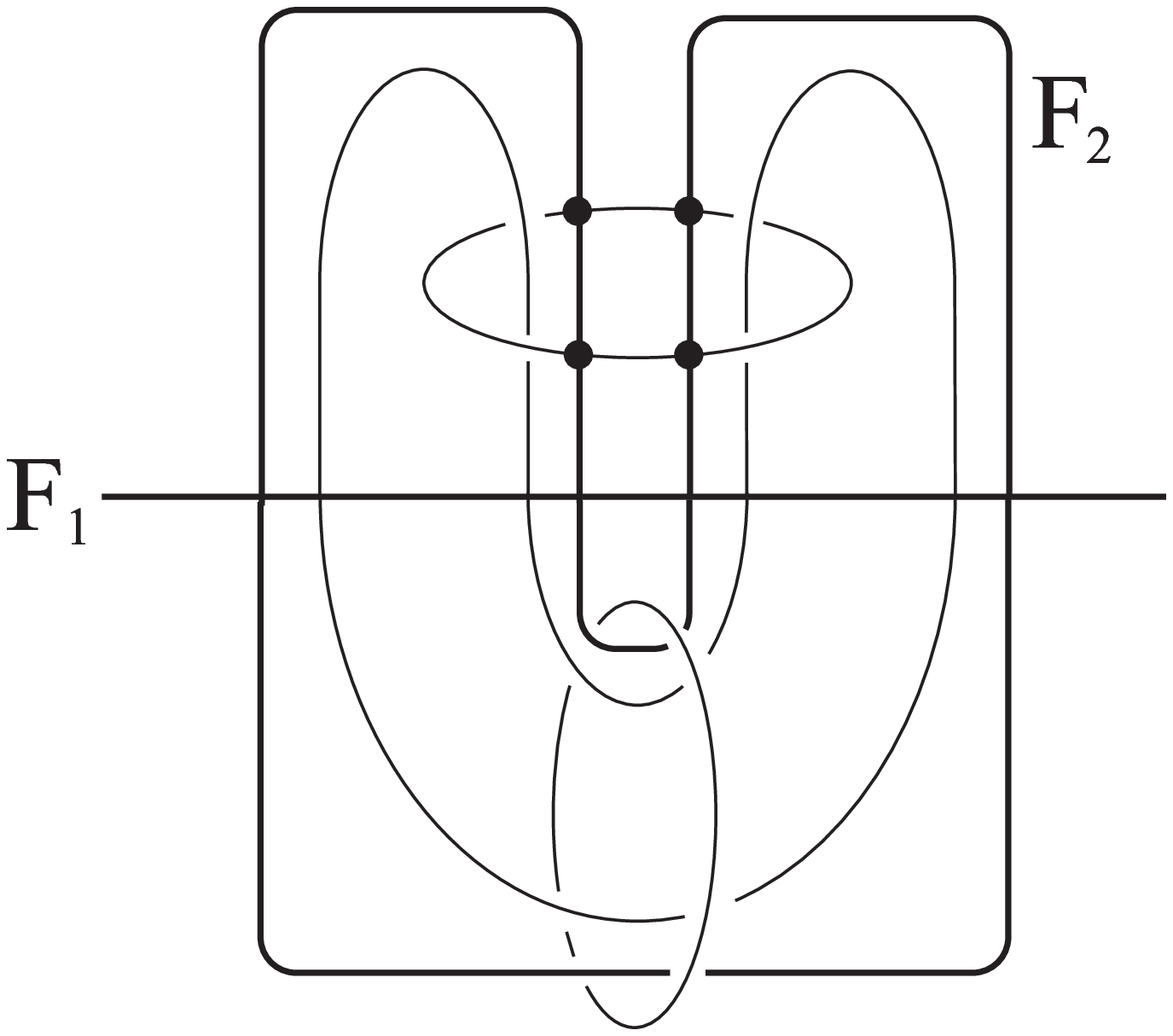} &
		\includegraphics[trim=0mm 0mm 0mm 0mm, width=.4\linewidth]{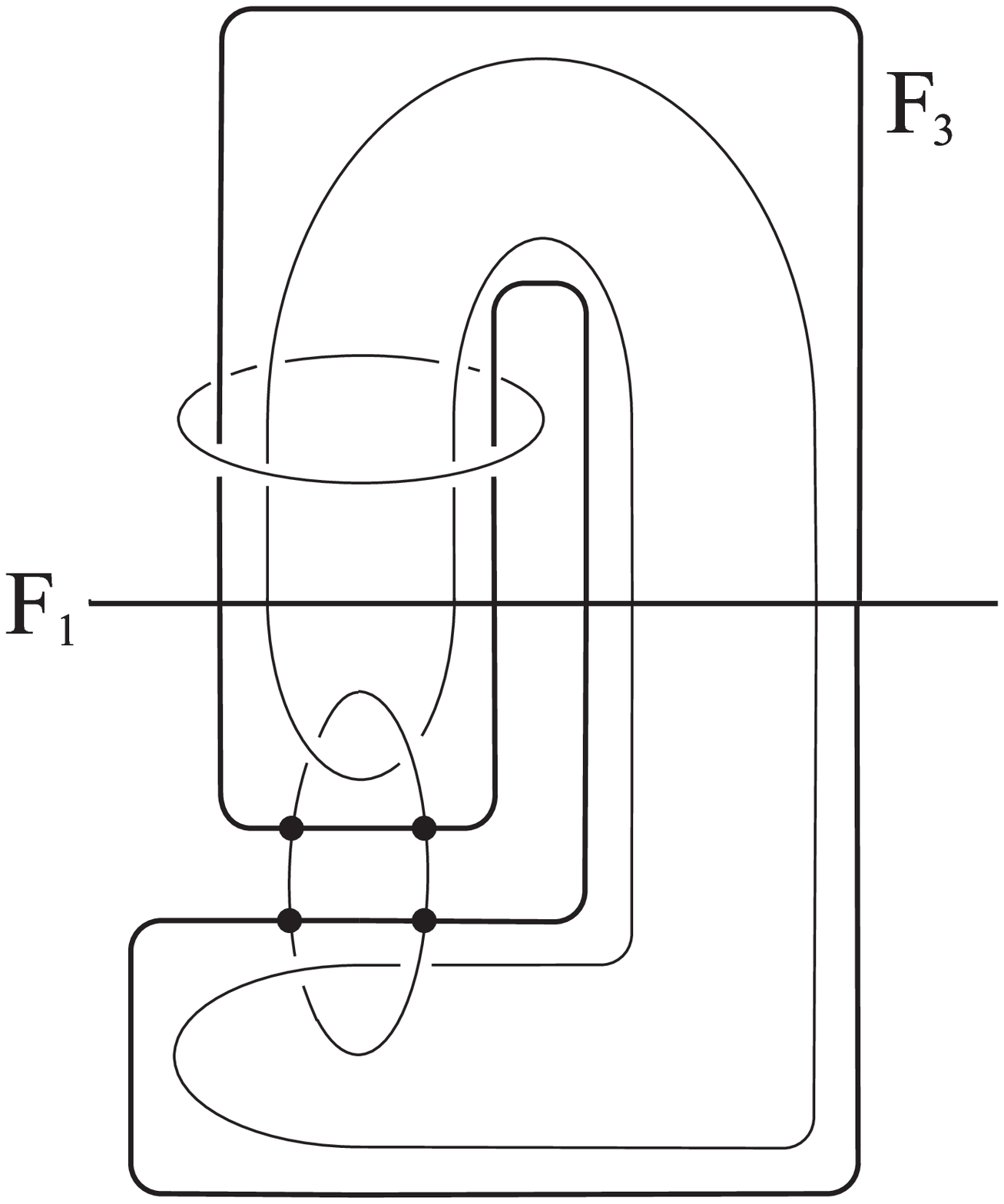}
	\end{tabular}
	\end{center}
	\caption{the Borromean rings and three Hopf tangle decomposing spheres}
	\label{Borromean_fig}
\end{figure}

\bigskip
\begin{example}\label{Borromean}
The Borromean rings is the only element of $\mathcal{L}(m_1,e_2)\cap \mathcal{L}(e_1,m_2)$ since it can be deformed as Figure \ref{Borromean_fig}.
By Remark \ref{class_remark}, it does not belong to $\mathcal{L}(m_1,m_2)$ and $\mathcal{L}(e_1,e_2)$.
Theorem \ref{Hopf} shows that the Borromean rings has exactly three Hopf tangle decomposing spheres $F_1$, $F_2$ and $F_3$, and that it is small.
Matsuda (\cite{HM}) also proved that the Borromean rings is small by using clasp disks bounded by the Borromean rings.
\end{example}

\section{Preliminary}

Let $S$ be a 2-sphere and $T$ a 1-manifold properly embedded in $S\times I$.
In this section, we consider incompressible surfaces embedded in $(S\times I,T)$, with the projection $p:S\times I\to S\times \{0\}$ and the height function $h:S\times I\to I$.
Let $F$ be a surface properly embedded in $S\times I$ such that $F$ is disjoint from $T$ or intersects $T$ transversely in the interior of $F$.
An isotopy $\phi_t$ of $S\times I$ is called a {\em horizontal isotopy} if $h\circ \phi_t=h$ for all $t\in [0,1]$.
A surface $F$ is said to be {\em horizontally $\partial$-parallel} {\em in} $S\times I$ if $F$ is $\partial$-parallel in $S\times I$ and there exists a horizontal isotopy $\phi_t$ of $S\times I$ such that $p\circ \phi_1$ is a homeomorphism on $F$.
We say that $F$ is {\em $T$-compressible} in $(S\times I,T)$ if there exists a disk $D$ in $S\times I$ such that $D\cap T=\partial D \cap T=\alpha$ is a subarc of $T$ and $D\cap F=\partial D -{\rm int} \alpha$.
A surface $F$ is {\em $T$-incompressible} in $(S\times I,T)$ if it is not $T$-compressible.

An $n$-dimensional submanifold $V$ in $S\times I$ is {\em vertical} if $p(V)$ is an $(n-1)$-dimensional manifold.
A submanifold $H$ in $S\times I$ is {\em horizontal} if $H$ is contained in a level 2-sphere $S\times \{t\}$ for some $t\in I$.

\bigskip
\begin{lemma}\label{horizontal}
{\it Let $S$ be a 2-sphere, $X$ a union of points in $S$, and 
$F\subset (S\times I,X\times I)$ an $(X\times I)$-incompressible surface.
If $F$ is horizontally $\partial$-parallel in $S\times I$, then $F$ is $\partial$-parallel in $(S\times I, X\times I)$.}
\end{lemma}

\bigskip
\begin{proof}
After a suitable horizontal isotopy of $S\times I$, we may assume that $p$ is a homeomorphism on $F$.
Moreover, by a suitable horizontal isotopy of $S\times I$ fixing $F$, we may assume that $p$ is a homeomorphism on $(X\times I)\cap V$ since $X\times I$ is monotone in $S\times I$, where $V$ is a 3-manifold bounded by $F$ and $p(F)$ in $S\times I$.
Then $p((X\times I)\cap V)$ consists of disjoint arcs embedded in $p(F)$.
See Figure \ref{parallel}.

\begin{figure}[htbp]
	\begin{center}
		\includegraphics[trim=0mm 0mm 0mm 0mm, width=.3\linewidth]{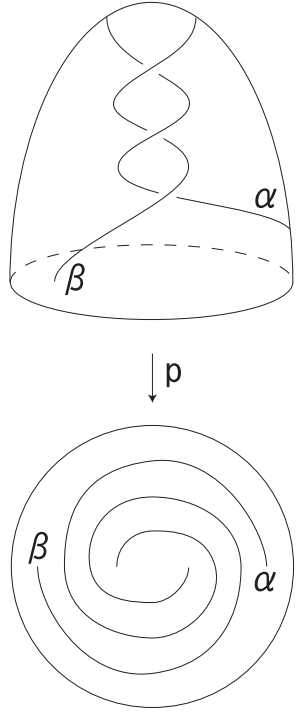}
	\end{center}
	\caption{$F$ is $\alpha$-compressible}
	\label{parallel}
\end{figure}

Suppose that there exists an arc $\alpha$ of $(X\times I)\cap V$ such that $\partial \alpha$ is contained in $F$.
We note that the vertical disk $p^{-1}(p(\alpha))\cap V$ is disjoint from $((X\times I)\cap V)-\alpha$.
Then, there exists a disk $D$ in $p^{-1}(p(\alpha))\cap V$ such that $D\cap \alpha=\partial D\cap \alpha=\alpha$ and $\partial D-{\rm int}\alpha\subset F$.
This shows that $F$ is an $(X\times I)$-compressible in $(S\times I,X\times I)$, a contradiction.
Hence, any arc of $(X\times I)\cap V$ connects $F$ and $p(F)$ monotonously, and $F$ is $\partial$-parallel in $(S\times I,X\times I)$.
\end{proof}


\bigskip
\begin{lemma}
\label{maximal}
{\it Let $S$ be a 2-sphere, $X$ a union of points in $S$, and $F\subset (S\times I,X\times I)$ an incompressible and meridionally incompressible surface in a Morse position with respect to $h$.
Suppose that the number of critical points of $F$ is minimal up to isotopy of $F$ in $(S\times I,X\times I)$.
If $F$ has a maximal critical point, then $F$ is $\partial$-parallel in $(S\times I,X\times I)$.}
\end{lemma}

\bigskip
\begin{proof}
If $F$ is horizontally $\partial$-parallel in $S\times I$, then by Lemma \ref{horizontal}, $F$ is $\partial$-parallel in $(S\times I,X\times I)$.
So, hereafter we assume that $F$ is not horizontally $\partial$-parallel in $S\times I$.

Let $a_0,\ldots,a_n$ be the critical points of $F$ which are labeled so that the corresponding critical values $t_i=h(a_i)$ satisfy $t_{i-1}<t_i$ for all $i$, and let $a_m$ be the lowest maximal critical point of $F$.
If $m=0$, then $F$ is a disk with only one maximal critical point $a_0$, and hence horizontally $\partial$-parallel in $S\times I$.
This contradicts the assumpution.
Hence, we have $m\ge 1$.

Let $F_{max}$ be the maximal horizontally $\partial$-parallel subsurface of $F$ containing $a_m$, that is, a component of $F\cap h^{-1}([t_{m-j}+\epsilon,t_m+\epsilon])$ which is horizontally $\partial$-parallel in $h^{-1}([t_{m-j}+\epsilon,t_m+\epsilon])$, for an integer $j\ (\le m)$ as large as possible and a fixed sufficiently small positive real number $\epsilon$.
See Figure \ref{maximal_fig}.

\begin{figure}[htbp]
	\begin{center}
		\includegraphics[trim=0mm 0mm 0mm 0mm, width=.7\linewidth]{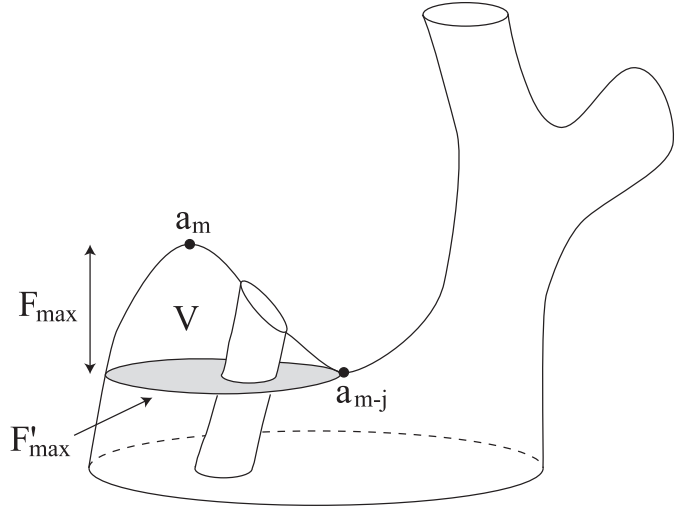}
	\end{center}
	\caption{maximal horizontally $\partial$-parallel subsurface $F_{max}$}
	\label{maximal_fig}
\end{figure}

Since $F_{max}$ is horizontally $\partial$-parallel in $S\times[t_{m-j}+\epsilon,t_m+\epsilon]$ and $(X\times[t_{m-j}+\epsilon,t_m+\epsilon])$-incompressible in $(S\times[t_{m-j}+\epsilon,t_m+\epsilon],X\times[t_{m-j}+\epsilon,t_m+\epsilon])$, by Lemma \ref{horizontal}, it is a $\partial$-parallel planar surface in $(S\times[t_{m-j}+\epsilon,t_m+\epsilon],X\times[t_{m-j}+\epsilon,t_m+\epsilon])$.
Therefore, by a suitable horizontal isotopy of $S\times I$, we may assume that $p$ is a homeomorphism on $F_{max}$.
Then, $F_{max}$ is isotopic to $F_{max}'=p^{-1}(p(F_{max}))\cap (S\times \{t_{m-j}+\epsilon\})$ in $(S\times [t_{m-j}+\epsilon,t_m+\epsilon], X\times [t_{m-j}+\epsilon,t_m+\epsilon])$, and $F_{max}$ and $F_{max}'$ bound a submanifold $V$ in $S\times [t_{m-j}+\epsilon,t_m+\epsilon]$ such that $p(F_{max})=p(F_{max}')=p(V)$.
We note that ${\rm int} V\cap F=\emptyset$ since $a_m$ is the lowest maximal critical point of $F$.
Moreover, by the minimality of the number of critical points of $F$, $F_{max}$ has only one maximal point $a_m$ and $(|\partial F_{max}|-1)$-saddle points.

Hereafter, we consider the critical point $a_{m-j}$.
If $a_{m-j}$ is a saddle point, then the saddle transformation at $a_{m-j}$ can be regarded as attaching a horizontal rectangular band $B$ to $F\cap (S\times [t_{m-j}+\epsilon, t_m+\epsilon])$ in a usual way.
Otherwise, since $a_m$ is the lowest maximal point, $a_{m-j}$ is a minimal point and can be regarded as attaching a horizontal disk $E$ to $F\cap (S\times [t_{m-j}+\epsilon, t_m+\epsilon])$.
Then, there are the following possibility for the critical point $a_{m-j}$.

\begin{description}
	\item [Case 1] The band $B$ is outside of $F_{max}'$ and connects a same boundary component of $F_{max}$.
	\item [Case 2] The band $B$ is outside of $F_{max}'$ and connects a boundary component of $F_{max}$ and a component of $F\cap (S\times [t_{m-j}+\epsilon,t_m+\epsilon])$ other than $F_{max}$.
	\item [Case 3] The band $B$ is inside of $F_{max}'$ and connects a same boundary component of $F_{max}$.
	\item [Case 4] The band $B$ is inside of $F_{max}'$ and connects two different boundary components of $F_{max}$.
	\item [Case 5] The disk $E$ caps $F_{max}$ along a boundary component of $F_{max}$.
\end{description}

In Case 1, $p$ is a homeomorphism on $F_{max}\cup B$.
Hence, a component of $F\cap h^{-1}([t_{m-j}-\epsilon,t_m+\epsilon])$ containing $a_m$ is horizontally $\partial$-parallel in $h^{-1}([t_{m-j}-\epsilon,t_m+\epsilon])$.
This contradicts the maximality of $F_{max}$.

In Case 2, by an isotopy of $F_{max}$ to $F_{max}'$ and pushing out it through $F_{max}'$, the maximal critical point $a_m$ and the saddle point $a_{m-j}$ are mutually canceled without increasing the number of critical points of $F$.
This contradicts the minimality of the number of critical points of $F$.

In Case 3, there exists a vertical compressing disk $D$ for $F_{max}\cup B$ in $(S\times [t_{m-j}-\epsilon, t_m+\epsilon], X\times [t_{m-j}-\epsilon, t_m+\epsilon])$ such that $D\cap (X\times I)=\emptyset$, $D\cap B=\partial D\cap B=\alpha$ is an arc through $B$, and $p(D)=p(\alpha)$.
By the incompressibility of $F$ in $(S\times I,X\times I)$, there exists a disk $D'$ in $F$ such that $\partial D'=\partial D$.
Then, by an isotopy of $D'$ to $D$ and a suitable isotopy of $F$, the saddle point $a_{m-j}$ is eliminated without increasing the number of critical points of $F$.
This contradicts the minimality of the number of critical points of $F$.

In Case 4, there exists a vertical compressing disk $D$ for $F_{max}\cup B$, as same as Case 3.
Actually, this disk $D$ gives a compressing disk for $F$ in $(S\times I,X\times I)$ since $\partial D$ is essential in $F$.
This contradicts the incompressibility of $F$ in $(S\times I,X\times I)$.

In Case 5, first we consider the case that $F_{max}$ is a disk.
If $p(F_{max})$ contains $p(E)$, then $F=F_{max}\cup E$ is a 2-sphere which bounds a 3-ball $V'$ in $S\times I$.
We note that $V'\cap (X\times I)=\emptyset$ since $E\cap (X\times I)=\emptyset$ and $F_{max}$ is isotopic to $E$ in $(V',(X\times I)\cap V')$.
Hence $F$ bounds a 3-ball $V'$ in $S\times I-X\times I$, and this contradicts the incompressibility of $F$ in $(S\times I,X\times I)$.
Otherwise, $F_{max}\cup E$ is isotopic to a level 2-sphere $h^{-1}(t_{m-j})$ in $(S\times I,X\times I)$, but this contradicts the assumpution that $F$ is not horizontally $\partial$-parallel in $S\times I$.

Next, we consider the case that $F_{max}$ is not a disk.
In this case, $p(F_{max})$ does not contain $p(E)$.
Hence, $F_{max}\cup E$ is isotopic to a planar surface in $h^{-1}(t_{m-j})$, and there exists an isotopy of $F_{max}\cup E$ so that it has only one maximal point $a_m$ and $(|\partial F_{max}|-2)$-saddle points since a saddle point of $F_{max}$ and a minimal point $a_{m-j}$ are mutually canceled.
This contradicts the minimality of the number of critical points of $F$.

Hence, any case does not occur, the lemma is proved.
\end{proof}

\bigskip
The next lemma is a special case of Lemma \ref{maximal}.
In the case that the number of strings is less than or equal to 5, incompressible and meridionally incompressible surfaces are very restricted.

\bigskip
\begin{lemma}
\label{5 points}
{\it Let $S$ be a 2-sphere, $X$ a union of points in $S$ less than or equal to 5, and $F\subset (S\times I,X\times I)$ an incompressible and meridionally incompressible surface in a Morse position with respect to $h$.
Suppose that the number of critical points of $F$ is minimal up to isotopy of $F$ in $(S\times I,X\times I)$.
Then, one of the following holds.
\begin{enumerate}
	\item $F$ has no critical point and any loop of $F\cap (S\times \{t\})$ does not bound a disk $D$ in $S\times \{t\}$ such that $|D\cap (X\times I)|\le 1$.
	\item $F$ is $\partial$-parallel in $(S\times I,X\times I)$.
\end{enumerate}}
\end{lemma}

\bigskip
\begin{proof}
If $F$ has a maximal or minimal critical point, then by Lemma \ref{maximal}, $F$ is $\partial$-parallel in $(S\times I, X\times I)$ and we have a conclusion 2 of Lemma \ref{5 points}.

Hereafter, we assume that $F$ has neither a maximal nor minimal critical point, and suppose that $F$ has a saddle point $a_s$.
Let $F_s$ be a component of $F\cap h^{-1}([h(a_s)-\epsilon, h(a_s)+\epsilon])$ containing $a_s$.
Then, at least one of boundary components of $F_s$ bounds a disk in $S\times \{h(a_s)\pm\epsilon\}$ which is disjoint from $X\times I$ or intersects $X\times I$ in one point.
Therefore, there exists a horizontal compressing disk or horizontal meridionally compressing disk $D$ for $F$ in $(S\times [h(a_s)-\epsilon, h(a_s)+\epsilon],X\times [h(a_s)-\epsilon, h(a_s)+\epsilon])$.
Since $F$ is incompressible and meridionally incompressible in $(S\times I, X\times I)$, there exists a disk $D'$ in $F$ such that $\partial D'=\partial D$ and $D\cup D'$ bounds a 3-ball in $S\times I-X\times I$ or a trivial $($3-ball, arc$)$-pair in $(S\times I,X\times I)$.
Moreover, since $D$ is horizontal, $D'$ contains at least one maximal or minimal critical point.
This contradicts the assumption and hence $F$ has no critical point.
Similarly, if there exists a loop of $F\cap (S\times \{t\})$ which bounds a disk $D$ in $S\times \{t\}$ such that $|D\cap (X\times I)|\le 1$, then it contradicts the assumption.
Therefore, we have a conclusion 1 of Lemma \ref{5 points}.
\end{proof}

\bigskip
\begin{lemma}\label{rational}
{\it Let $(B,T)$ be a rational tangle, $F$ a surface properly embedded in $B$ which is disjoint from $T$ or intersects $T$ transversely.
If $F$ is incompressible and meridionally incompressible in $(B,T)$, then one of the following holds.
\begin{enumerate}
	\item $F$ is $\partial$-parallel in $(B,T)$.
	\item $F$ is a disk separating two strings of $T$.
\end{enumerate}}
\end{lemma}

\bigskip
\begin{proof}
Let $B$ be the north hemisphere $\{(x,y,z,w)\in \Bbb{R}^4|x^2+y^2+z^2+w^2=1, w\ge 0\}$, and $h:B\to \Bbb{R}$ be the height function defined by $(x,y,z,w)\mapsto w$.
Since $(B,T)$ is a rational tangle, $T$ can be isotoped so that the restriction of $h$ to $T$ is a Morse function with two maximal critical points.
Let $S$ a level 2-sphere below two maximal points of $T$, and put $M_{thin}=B\cap h^{-1}([h(S),1])$ and $M_{thick}=B\cap h^{-1}([0,h(S)])$.
See Figure \ref{rational_fig}.

\begin{figure}[htbp]
	\begin{center}
		\includegraphics[trim=0mm 0mm 0mm 0mm, width=.5\linewidth]{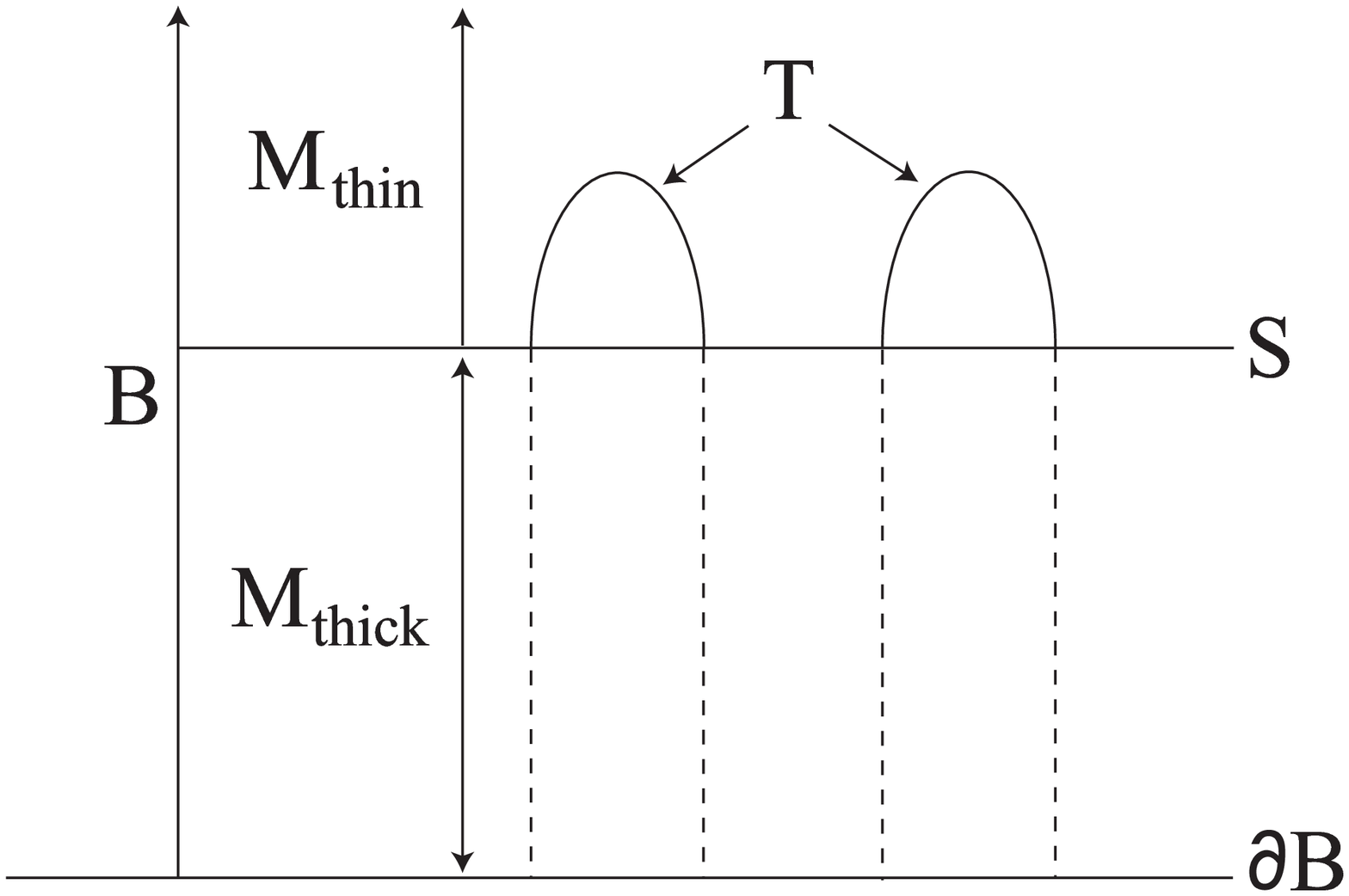}
	\end{center}
	\caption{$M_{thin}$-$M_{thick}$ decomposition of $(B,T)$}
	\label{rational_fig}
\end{figure}

We isotope $F$ so that $F$ is a Morse position with respect to $h$.
In the same way as Claim \ref{thin region}, we may assume that $F\cap M_{thin}$ consists of incompressible disks which separate two strings of $T\cap M_{thin}$.
Hereafter, we suppose that $|F\cap M_{thin}|$ and the critical points of $F$ are minimal up to isotopy of $F$ in $(B,T)$.
Then each component of $F\cap M_{thin}$ is an incompressible disk with only one maximal critical point, and by the proof of Claim \ref{thick region}, each component of $F\cap M_{thick}$ is incompressible and meridionally incompressible in $(M_{thick}, T\cap M_{thick})$.
Hence, by Lemma \ref{5 points}, $F\cap M_{thick}$ has no critical point or $F$ is $\partial$-parallel in $(M_{thick}, T\cap M_{thick})$.

In the latter case, if $F\cap M_{thick}$ is $\partial$-parallel into $S$, then  $F$ can be isotoped so that it is entirely contained in $M_{thin}$.
Since $(M_{thin},T\cap M_{thin})$ is a trivial tangle, $F$ is compressible or meridionally compressible in $(M_{thin},T\cap M_{thin})$.
This contradicts the incompressibility or meridionally incompressibility of $F$ in $(B,T\cap B)$.
Otherwise, $F\cap M_{thick}$ is $\partial$-parallel into $\partial B$, and we obtain the conclusion 1 of Lemma \ref{rational}.

In the formar case, $F\cap M_{thick}$ has no critical point and for any level 2-sphere $P$ in $M_{thick}$, any loop of $F\cap P$ separates four points of $T\cap P$ into two points and two points in $P$.
It follows that $F\cap M_{thick}$ is a vertical annulus disjoint from $T$.
Hence, $F$ is a disk separating two strings of $T$, and we have the conclusion 2 of Lemma \ref{rational}.
%
\end{proof}

\section{Proof}



\bigskip
\begin{proof} (of Theorem \ref{main})
First, we focus the thin regions $M_{thin}=M_1\cup M_2\cup \cdots \cup M_{k-1}\cup M_k$, where $M_1$ is the top 3-ball, $M_2,\ldots,M_{k-1}$ are thin regions which contain thin level spheres $S_2,\ldots,S_{k-1}$ respectively, and $M_k$ is the bottom 3-ball, where $h(S_i)>h(S_{i+1})$ for $i=2,\ldots,k-2$.
Put $S_1=h^{-1}(1-\epsilon)$ and $S_k=h^{-1}(-1+\epsilon)$ for a sufficiently small positive real number $\epsilon$.
Then, for each $i$, $S_i$ separates $M_i$ into two submanifolds $M_i^+$ and $M_i^-$, where $M_i^+$ is the upper part and $M_i^-$ is the lower part.
Put $K\cap M_i^{\pm}=\alpha_i^{\pm}\cup \beta_i^{\pm}$, where $\alpha_i^{\pm}$ is a union of arcs each of which has only one critical point, and $\beta_i^{\pm}$ is a union of vertical arcs.
See Figure \ref{thin1}.

\begin{figure}[htbp]
	\begin{center}
		\includegraphics[trim=0mm 0mm 0mm 0mm, width=.7\linewidth]{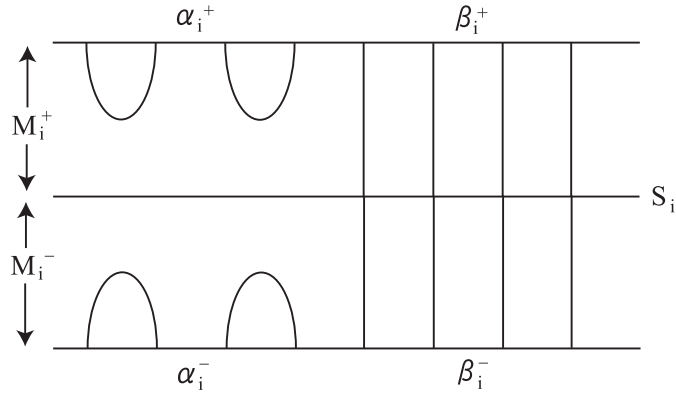}
	\end{center}
	\caption{$M_i^{\pm}$, $S_i$, $\alpha_i^{\pm}$, $\beta_i^{\pm}$}
	\label{thin1}
\end{figure}

\bigskip
\begin{claim}
\label{thin region}
By an isotopy of $F$ in $(S^3,K)$, we may assume that each component of $F\cap M_i^{\pm}$ is either an incompressible disk in $(M_i^{\pm},K\cap M_i^{\pm})$ with only one critical point or an incompressible vertical annulus in $(M_i^{\pm},K\cap M_i^{\pm})$ disjoint from $K$.
Moreover, if $F$ is meridionally incompressible in $(S^3,K)$, then we may assume that each component of $F\cap M_i^{\pm}$ is meridionally incompressible in $(M_i^{\pm},K\cap M_i^{\pm})$.
\end{claim}

\bigskip
\begin{proof} (of Claim \ref{thin region})
By an incompressibility of $F$ in $(S^3,K)$, we may assume that $F\cap S_i$ consists of essential loops in $S_i-K$.
Let $T_i^{\pm}\subset M_i^{\pm}$ be a union of mutually disjoint monotone arcs connecting a minimal/maximal critical point of $K\cap M_i^{\pm}$ to a point in $S_i-(F\cap S_i)$.
We may assume that $F$ and $T_i^{\pm}$ are in general position.
Put $N_i^{\pm}=N(S_i\cup T_i^{\pm};M_i^{\pm})$.
Then, each component of $F\cap N_i^{\pm}$ is either an incompressible disk in $(N_i^{\pm},K\cap N_i^{\pm})$ or an incompressible annulus in $(N_i^{\pm},K\cap N_i^{\pm})$ disjoint from $K$.
Since $(M_i^{\pm}-{\rm int} N_i^{\pm},K\cap (M_i^{\pm}-{\rm int} N_i^{\pm}))$ admits a product structure,
there exists an isotopy of $\partial N_i^{\pm}-S_i$ into $\partial M_i^{\pm}-S_i$ in $(M_i^{\pm},K\cap M_i^{\pm})$ keeping such conditions.
Hence, we have that each component of $F\cap M_i^{\pm}$ is either an incompressible disk in $(M_i^{\pm},K\cap M_i^{\pm})$ with only one critical point or an incompressible vertical annulus in $(M_i^{\pm},K\cap M_i^{\pm})$ disjoint from $K$.
See Figure \ref{thin}.

Moreover, suppose that $F$ is meridionally incompressible in $(S^3,K)$.
If there exists a component of $F\cap M_i^{\pm}$ which is meridionally compressible in $(M_i^{\pm},K\cap M_i^{\pm})$, then it is a vertical annulus.
We choose an innermost vertical annlus $A$ and a meridionally compressing disk $D$ for $A$ in $(M_i^{\pm},K\cap M_i^{\pm})$.
Then, there exists a disk $D'$ in $F$ such that $\partial D'=\partial D$ and $D\cup D'$ bounds a trivial $($3-ball, arc$)$-pair in $(S^3,K)$.
An isotopy of $F$ bringing $D'$ to $D$ in $(S^3,K)$ turns the vertical annulus $A$ to a meridionally compressing disk $A'$ in $(M_i^{\pm},K\cap M_i^{\pm})$, and it can be pushed out from $M_i^{\pm}$.
Here, if $A'$ enters $M_i^{\mp}$, then it can be pushed out from $M_i$ entirely.
Hence, each component of $F\cap M_i^{\pm}$ is meridionally incompressible in $(M_i^{\pm},K\cap M_i^{\pm})$.
\end{proof}

\begin{figure}[htbp]
	\begin{center}
		\includegraphics[trim=0mm 0mm 0mm 0mm, width=.8\linewidth]{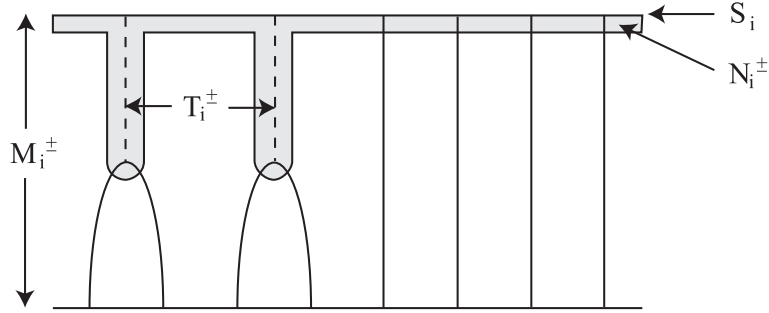}
	\end{center}
	\caption{$M_i^{\pm}$, $S_i$, $T_i^{\pm}$ and $N_i^{\pm}$}
	\label{thin}
\end{figure}

Hence, $F\cap M_i$ consists of incompressible disks with only one minimal/maximal critical point or vertical incompressible annuli in $(M_i,K\cap M_i)$.
Moreover, if $F$ is meridionally incompressible in $(S^3,K)$, then $F\cap M_i$ is also meridionally incompressible in $(M_i,K\cap M_i)$.
Hereafter, we suppose that $|F\cap M_{thin}|$ is minimal up to isotopy of $F$ in $(S^3,K)$ under such conditions, and moreover that the number of critical points of $F\cap M_{thick}$ is minimal among all Morse positions of $F$.
We note that $|F\cap M_{thin}|\ne 0$ since if $F\subset M_{thick}$, then by Lemma \ref{maximal}, $F$ is parallel to a thick level sphere and it is compressible.

\bigskip
\begin{claim}\label{thick region}
Each component of $F\cap M_{thick}$ is incompressible and meridionally incompressible in $(M_{thick},K\cap M_{thick})$.
\end{claim}

\bigskip
\begin{proof} (of Claim \ref{thick region})
Suppose that there exists a component of $F\cap M_{thick}$ which is compressible in $(M_{thick},K\cap M_{thick})$.
Then, there exists a compressing disk $D$ for $F\cap M_{thick}$ in $(M_{thick},K\cap M_{thick})$.
Since $F$ is incompressible in $(S^3,K)$, there exists a disk $D'$ in $F$ such that $\partial D'=\partial D$.
Then, an isotopy of $F$ from $D'$ to $D$ decreases $|F\cap M_{thin}|$ since $D'$ contains at least one component of $F\cap M_{thin}$.
This is a contradiction to the minimality of $|F\cap M_{thin}|$.
Finally, in the case that $F$ is meridionally incompressible in $(S^3,K)$, the proof is similar to the case that $F$ is incompressible.
\end{proof}

\bigskip
By the above, the conditions (1), (2) and (4) in Section 1 are satisfied.
Moreover, if $F\cap M_{thick}$ does not contain a maximal or minimal critical point, then the condition (3) is satisfied.
Hence $F$ is in an essential Morse position related to $K$, and we have the conclusion 2 of Theorem \ref{main}.

Hereafter, we suppose that $F\cap M_{thick}$ has a maximal or minimal critical point, and we will show that $F$ is isotopic to a thin level sphere, which is the conclusion 1 of Theorem \ref{main}.
Without loss of generality, let $a_m$ be the lowest maximal critical point of $F\cap M_{thick}$ in a thick region.
Then, by Lemma \ref{maximal}, a component $F_{max}$ of $F\cap M_{thick}$ containing $a_m$ is $\partial$-parallel in $(M_{thick}, K\cap M_{thick})$.
If $F_{max}$ is parallel into $\partial M_k$, then $F$ is isotoped so that it is entirely contained in the trivial tangle $(M_k,K\cap M_k)$, and it is compressible in the trivial tangle.
Therefore, without loss of generality, we may assume that $F_{max}$ is isotopic to a planar surface contained in $\partial M_i^+-S_i$ for a thin region $M_i$.
Put $F\cap M_i^{\pm}=D_i^{\pm}\cup A_i^{\pm}$, where $D_i^{\pm}$ is a union of incompressible disks with only one critical point and $A_i^{\pm}$ is a union of vertical annuli in $(M_i^{\pm},K\cap M_i^{\pm})$.
Since each arc of $\alpha_i^{\pm}$ has only one critical point, there exists a union $\Delta_i^{\pm}$ of vertical disks in $M_i^{\pm}$ for $\alpha_i^{\pm}$ such that $\Delta_i^{\pm}\cap \alpha_i^{\pm}=\partial \Delta_i^{\pm}\cap \alpha_i^{\pm}=\alpha_i^{\pm}$ and $\gamma_i^{\pm}=\partial \Delta_i^{\pm} - {\rm int} \alpha_i^{\pm}\subset \partial M_i^{\pm}-S_i$.
Moreover, by a horizontal isotopy, we can take $\Delta_i^{\pm}$ so that $\Delta_i^{\pm}\cap (D_i^{\pm}\cup A_i^{\pm})=\emptyset$ since each component of $D_i^{\pm}$ has only one critical point and $A_i^{\pm}$ is vertical.
See Figure \ref{thin2}.

\begin{figure}[htbp]
	\begin{center}
		\includegraphics[trim=0mm 0mm 0mm 0mm, width=.7\linewidth]{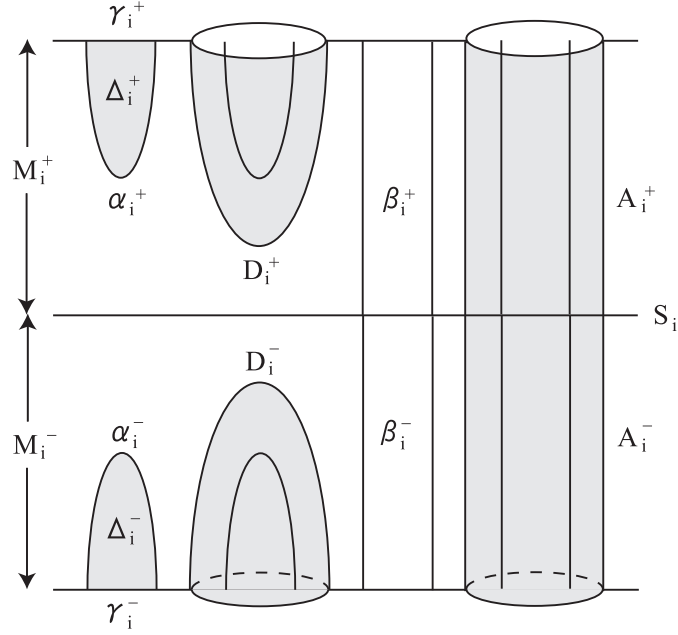}
	\end{center}
	\caption{$M_i^{\pm}$, $S_i$, $\alpha_i^{\pm}$, $\beta_i^{\pm}$, $\gamma_i^{\pm}$, $D_i^{\pm}$ and $A_i^{\pm}$}
	\label{thin2}
\end{figure}

Now, we isotope $F_{max}$ so that it is contained in $\partial M_i^+-S_i$.
If $F_{max}$ contains a component of $\gamma_i^+$, then $F_{max}$ is $\alpha_i^+$-compressible and hence $F$ is $K$-compressible.
This contradicts the incompressibility of $F$ in $(S^3,K)$ and we have $F_{max}\cap \Delta_i^+=\emptyset$.
Let $F_{max}'$ be a component of $F_{max}\cup D_i^+\cup A_i^+\cup A_i^-$ containing $F_{max}$ and push it into $M_i$ slightly.
Since $\Delta_i^{\pm}\cap (F_{max}\cup D_i^{\pm}\cup A_i^{\pm})=\emptyset$, we may assume that $F_{max}'$ is contained in $M_i'=M_i-{\rm int} N(\Delta_i^{\pm};M_i)$.
We note that $(M_i',K\cap M_i')$ admits a product structure $(M_i',\beta_i^+\cup \beta_i^-)$.

There are following two cases for $F_{max}'$.

\begin{description}
	\item[Case 1] $\partial F_{max}\cap A_i^+=\emptyset$
	\item[Case 2] $\partial F_{max}\cap A_i^+\ne\emptyset$
\end{description}

In Case 1, $F$ is a 2-sphere which is entirely contained in $(M_i',\beta_i^+\cup \beta_i^-)$.
By Lemma \ref{maximal}, $F$ is $\partial$-parallel in $(M_i',\beta_i^+\cup \beta_i^-)$.
It follows that $F$ is isotopic to a thin level sphere $S_i$, and we have the conclusion 1 of Theorem \ref{main}.

In Case 2, $F_{max}'$ is a horizontally $\partial$-parallel planar surface in $(M_i',\beta_i^+\cup \beta_i^-)$.
By Lemma \ref{horizontal}, $F_{max}'$ is $\partial$-parallel in $(M_i',\beta_i^+\cup \beta_i^-)$.
We isotope $F_{max}'$ in $(M_i',\beta_i^+\cup \beta_i^-)$ so that it is contained in $M_i^-$.
Then, there exists a union $\delta$ of vertical $\partial$-compressing disks for $F_{max}'$ in $(M_i^-,K\cap M_i^-)$ such that $\partial \delta$ cuts $F_{max}'$ into a single disk, and $\delta\cap (D_i^-\cup \Delta_i^-)=\emptyset$.
By a $\partial$-compression of $F_{max}'$ along $\delta$ and a suitable isotopy, $F_{max}'$ turns into a disk with only one maximal critical point in $M_i^-$.
During the transformation above, $|F\cap M_i|$ does not increase, but at least one maximal critical point $a_m$ of $F\cap M_{thick}$ is eliminated.
This contradicts the supposition of the number of critical points of $F\cap M_{thick}$.
Hence, Theorem \ref{main} is proved.
\end{proof}

\bigskip
\begin{proof} (of Theorem \ref{Montesinos})
Let $(B,T)$ be a Montesinos tangle $T(r_1,r_2)$ with $r_i\ne\infty$ for $i=1,2$, and $D$ a disk dividing $(B,T)$ into two rational tangles $(B_1,T_1)$ and $(B_2,T_2)$ so that $D=B_1\cap B_2=\partial B_1\cap \partial B_2$ and $(B_i,T_i)$ has a slope $r_i$.
Let $F$ be a closed incompressible surface properly embedded in $B-{\rm int} N(T)$.
If $F$ is meridionally compressible in $(B,T)$, then we perform meridionally compressions as possible and obtain a disjoint union $F'$ of closed incompressible and meridionally incompressible surfaces.
Conversely, $F$ is obtained from $F'$ by meridional tubings along $T$.
Note that $F'\subset {\rm int} B$ since $F$ is a closed surface properly embedded in $B-{\rm int} N(T)$.

\bigskip
\begin{claim}\label{Montesinos1}
Let $S$ be a closed incompressible and meridionally incompressible surface in $(B,T)$ disjoint from $T$ or intersecting $T$ transversely.
Then $S$ is $\partial$-parallel in $(B,T)$.
\end{claim}

\bigskip
\begin{proof} (of Claim \ref{Montesinos1})
Suppose that $|S\cap D|$ is minimal up to isotopy of $S$.
Then, by the incompressibility and meridionally incompressibility of $S$ in $(B,T)$, each component of $S\cap B_i$ is incompressible and meridionally incompressible in $(B_i,T_i)$.
We note that $S\cap D\ne \emptyset$ since a rational tangle contains no closed incompressible and meridionally incompressible surface, and that any loop of $S\cap D$ is parallel to $\partial D$ in $D-T$.
By Lemma \ref{rational}, each component $P$ of $S\cap B_i$ is $\partial$-parallel in $(B_i,T_i)$ or a disk separating two strings of $T_i$.
If there exists a component $P$ of $S\cap B_i$ of the latter type, then $r_i=\infty$ and this contradicts the hypothesis for $r_i$.
Hence, by the minimality of $|S\cap D|$, both of $S\cap B_1$ and $S\cap B_2$ consist of $\partial$-parallel disks that are not parallel into $D$.
By connecting these components of $S\cap B_i$, we have a 2-sphere $S$ which is $\partial$-parallel in $(B,T)$.
\end{proof}

By Claim \ref{Montesinos1}, $F'$ is a $\partial$-parallel 2-sphere in $(B,T)$.
Now, we recover $F$ from $F'$ by tubing $F'$ along $T$.
Then $F$ is $\partial$-parallel in $B-{\rm int} N(T)$, hence the Montesinos tangle $T(r_1,r_2)$ is small.
\end{proof}

\bigskip
\begin{remark}\label{length}
In general, Montesinos tangles are not small.
Tsutsumi informed me that there are Montesinos tangles with length more than two that are not small.
Actually, any closed incompressible and meridionally incompressible surface in a Montesinos tangle $T(r_1,\ldots,r_n)$ are parallel to $\partial (B_i\cup \cdots\cup B_j)$, where $B_i,\ldots,B_j$ are successive two or above rational tangles.
Then, we perform a meridional tubing along $T_i\cup\cdots\cup T_j$ and obtain a closed incompressible surface which is not $\partial$-parallel in the exterior of $T(r_1,\ldots,r_n)$.
\end{remark}

\bigskip
\begin{proof} (of Theorem \ref{theta})
Let $G$ be a 2-bridge theta curve or handcuff graph, $F$ an incompressible and meridionally incompressible closed surface in $(S^3,G)$.
In the same way as the case of knots, we decompose $S^3$ into $M_{thin}=M_1\cup M_2$ and $M_{thick}$, where $M_1$ is the top 3-ball and $M_2$ is the bottom 3-ball.
Then, one vertex and the maximal critical point of $G$ are contained in $M_1$, and another vertex and the minimal critical point of $G$ are contained in $M_2$.
Moreover, by the proof of Theorem \ref{main}, we may assume that $F\cap M_1$ are consists of incompressible disks in $(M_1,G\cap M_1)$ each of which has only one maximal point and disjoint from $G$, $F\cap M_{thick}$ consists of incompressible and meridionally incompressible surfaces each of which has no minimal/maximal critical point, $F\cap M_2$ are consists of incompressible disks in $(M_2,G\cap M_2)$ each of which has only one minimal point and disjoint from $G$.
Suppose that $|F\cap M_{thin}|$ and the number of critical points of $F\cap M_{thick}$ is minimal up to isotopy of $F$ in $(S^3,G)$.
By Lemma \ref{5 points}, $F\cap M_{thick}$ has no critical point or $F$ is $\partial$-parallel in $(M_{thick},G\cap M_{thick})$.

In the latter case, $F$ can be isotoped so that it is entirely contained in $M_{thin}$.
Then, $F$ is compressible or meridionally compressible in $(M_{thin},G\cap M_{thin})$ since $G\cap M_{thin}$ is contained in a disk properly embedded in $M_{thin}$.
This is a contradiction.

In the formar case, $F\cap M_{thick}$ has no critical point and for any thick level 2-sphere $P$ for $G$, any loop of $F\cap P$ separates five points of $G\cap P$ into two points and three points in $P$.
It follows that $F\cap M_{thick}$ is a vertical annulus, and $F$ is a sphere with one maximum and one minimum.
Hence, $F$ is a 2-sphere illustrated in Figure \ref{theta-figure}.
\end{proof}

\bigskip
\begin{proof} (of Theorem \ref{Hopf})
Let $L$ be a link in $S^3$ which admits a Hopf tangle decomposition $(S^3,L)=(B_1,T_1)\cup (B_2,T_2)$, and $F$ be an incompressible and meridionally incompressible surface in $(S^3,L)$.
For each $i=1,2$, the trivial loop $C_i\subset T_i$ bounds a disk $D_i$ in $B_i$ which intersects each string of $T_i$ in only one point.
Suppose that $|F\cap (D_1\cup D_2)|$ is minimal up to isotopy of $F$.
By the incompressibility and meridional incompressibility of $F$, each component of $F\cap D_i$ is either an arc separating two points of $T_i\cap D_i$ in $D_i$ or a loop parallel to $\partial D_i$ in $D_i-T_i$.
Put $H_i=N(D_i)$.
Then, $(H_i,T_i\cap H_i)$ is also a Hopf tangle and each component of $F\cap H_i$ is either a ``meridian'' disk bounded by the meridian of $(H_i,T_i\cap H_i)$ or a ``vertical'' annulus whose boundary component is parallel to the equator of $(H_i,T_i\cap H_i)$ in $\partial H_i-\partial (T_i\cap H_i)$.
See Figure \ref{component}.
We remark about $F\cap H_i$ that a meridian disk is not compatible with a vertical annulus.

\begin{figure}[htbp]
	\begin{center}
		\includegraphics[trim=0mm 0mm 0mm 0mm, width=.7\linewidth]{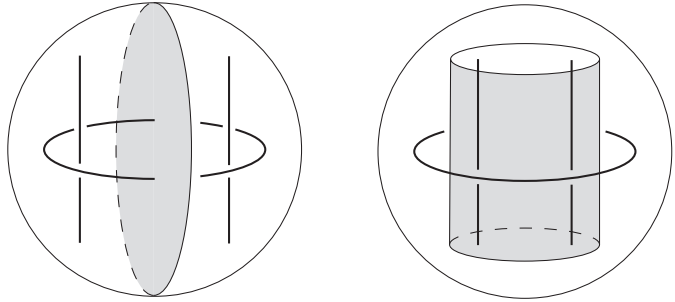}
	\end{center}
	\caption{a ``meridian'' disk and a ``vertical'' annulus}
	\label{component}
\end{figure}

The remainder of these two small Hopf tangles $S^3-\rm{int}(H_1\cup H_2)$ admits a product structure $(S\times I, X\times I)$, where $S$ is the decomposing sphere for $(S^3,L)$ and $X$ is four points of $L\cap S$.
Under the minimal condition of $|F\cap (D_1\cup D_2)|$, we put $F\cap (S\times I)$ in a Morse position with respect to the height function $h:S\times I\to I$ and suppose that the number of critical points of $F\cap (S\times I)$ is minimal among all Morse positions of $F\cap (S\times I)$.
In the same way as Claim \ref{thick region}, each component of $F\cap (S\times I)$ is incompressible and meridionally incompressible in $(S\times I,X\times I)$.
It follows from Lemma \ref{5 points} that each component of $F\cap (S\times I)$ is either a vertical annulus in $(S\times I,X\times I)$ which separates four strings of $X\times I$ into two and two, or $\partial$-parallel in $(S\times I,X\times I)$.
In the latter case, $F$ is parallel to the decomposing sphere $S$ when $F\cap (D_1\cup D_2)=\emptyset$, or it contradicts the minimality of $|F\cap (D_1\cup D_2)|$ or the incompressibility and meridional incompressibility of $F\cap (S\times I)$.
In the former case, $L$ belongs to $\mathcal{L}(m_1,m_2)\cup \mathcal{L}(m_1,e_2)\cup \mathcal{L}(e_1,m_2)\cup \mathcal{L}(e_1,e_2)$ since any loop of $F\cap \partial H_i$ is either the meridian or the equator of $(H_i,T_i\cap H_i)$ and each component of $F\cap (S\times I)$ is a vertical annulus disjoint from $X\times I$.
Thus the conclusion 1 of Theorem \ref{Hopf} holds.

Next, if $L\in \mathcal{L}(m_1,m_2)\cup \mathcal{L}(m_1,e_2)\cup \mathcal{L}(e_1,m_2)$, then there exists an another Hopf tangle decomposing sphere which is made of meridian disks and vertical annuli in $H_i$ and vertical annuli in $S\times I$ as Figure \ref{class}.
Coversely, if there exists a Hopf tangle decomposing sphere different from $S$, then it is an incompressible and meridionally incompressible closed surface in $(S^3,L)$.
Hence, by the conclusion 1 of Theorem \ref{Hopf}, we have $L\in \mathcal{L}(m_1,m_2)\cup \mathcal{L}(m_1,e_2)\cup \mathcal{L}(e_1,m_2)\cup \mathcal{L}(e_1,e_2)$.
But for a link which belongs to $\mathcal{L}(e_1,e_2)-\mathcal{L}(m_1,m_2)$, there does not exist an another Hopf tangle decomposing sphere, therefore we have $L\in \mathcal{L}(m_1,m_2)\cup \mathcal{L}(m_1,e_2)\cup \mathcal{L}(e_1,m_2)$.
Thus the first sentence of conclusion 2 of Theorem \ref{Hopf} holds.
Moreover, since any Hopf tangle decomposing sphere different from $S$ is made of meridian disks and vertical annuli in $H_i$ and vertical annuli in $S\times I$ as Figure \ref{class}, $L\in \mathcal{L}(m_1,m_2)\cup \mathcal{L}(m_1,e_2)\cup\mathcal{L}(e_1,m_2)$ has exactly two Hopf tangle decompositions, except the Borromean rings has exactly three Hopf tangle decompositions as Example \ref{Borromean}.
Thus the conclusion 2 of Theorem \ref{Hopf} holds.

Finally, suppose that $L$ is not small and let $Q$ be an incompressible and not $\partial$-parallel closed surface in $E(L)$.
Then, by meridional compressing $Q$ as possible, we get an incompressible and meridionally incompressible closed surface $Q'$ in $(S^3,L)$.
By the conclusion 1 of Theorem \ref{Hopf}, we have $L\in \mathcal{L}(m_1,m_2)\cup \mathcal{L}(m_1,e_2)\cup \mathcal{L}(e_1,m_2)\cup \mathcal{L}(e_1,e_2)$.
But if $L\not\in \mathcal{L}(e_1,e_2)$, then any meridional tubing $Q'$ cannot be performed since $Q'$ must be a Hopf tangle decomposing sphere.
Hence, $L\in \mathcal{L}(e_1,e_2)$ and $Q$ is an essential torus illustrated in Figure \ref{class}.
Thus the conclusion 3 of Theorem \ref{Hopf} holds.
\end{proof}

\bigskip
The following exercise can be proved by the similar method as Lemma \ref{rational}. So, we omit the proof.

\bigskip
\begin{xca}
Let $(B,T)$ be a Hopf tangle and $F$ a surface properly embedded in $B$ which is disjoint from $T$ or intersects $T$ transversely.
If $F$ is incompressible and meridionally incompressible in $(B,T)$, then one of the following holds.
\begin{enumerate}
	\item $F$ is a meridian disk for $(B,T)$.
	\item $F$ is a vertical annulus for $(B,T)$.
	\item $F$ is $\partial$-parallel in $(B,T)$.
\end{enumerate}

\end{xca}

\noindent{\bf Acknowledgement.}
The author would like to thank Yukihiro Tsutsumi for informing Remark \ref{length}, and the referee for pointing out the third Hopf tangle decomposing sphere for the Borromean rings.

\bibliographystyle{amsplain}

\end{document}